\documentclass[a4paper,12pt]{article}
\usepackage[all]{xy}
\usepackage[T1]{fontenc}
\usepackage[dvips]{graphicx}
\usepackage{amsmath,amssymb,amsfonts,amsthm,latexsym,mathrsfs,textcomp,verbatim, enumerate, soul, color, MnSymbol}
\xyoption{web}

\title{Priestley-type dualities for\\ partially ordered structures}
\author{{Olivia Caramello} \vspace{3 mm}\\ {\small DPMMS, University of Cambridge,}\\{\small Wilberforce Road, Cambridge CB3 0WB, UK}\\{\small O.Caramello@dpmms.cam.ac.uk}}
\date{March 13, 2012}
\begin{document}
\bgroup           
\let\footnoterule\relax  
\maketitle
\begin{abstract}
We introduce a general framework for generating dualities between categories of partial orders and categories of ordered Stone spaces; we recover in particular the classical Priestley duality for distributive lattices and establish several other dualities for different kinds of partially ordered structures. 
\end{abstract} 
\egroup 
\vspace{5 mm}

\tableofcontents 


\def\Monthnameof#1{\ifcase#1\or
   January\or February\or March\or April\or May\or June\or
   July\or August\or September\or October\or November\or December\fi}
\def\today{\number\day~\Monthnameof\month~\number\year}

%
%
%
\def\pushright#1{{
   \parfillskip=0pt            
   \widowpenalty=10000         
   \displaywidowpenalty=10000  
   \finalhyphendemerits=0      
  %
   \leavevmode                 
   \unskip                     
   \nobreak                    
   \hfil                       
   \penalty50                  
   \hskip.2em                  
   \null                       
   \hfill                      
   {#1}                        
  %
   \par}}                      

\def\qed{\pushright{$\square$}\penalty-700 \smallskip}

\newtheorem{theorem}{Theorem}[section]

\newtheorem{proposition}[theorem]{Proposition}

\newtheorem{scholium}[theorem]{Scholium}

\newtheorem{lemma}[theorem]{Lemma}

\newtheorem{corollary}[theorem]{Corollary}

\newtheorem{conjecture}[theorem]{Conjecture}

\newenvironment{proofs}%
 {\begin{trivlist}\item[]{\bf Proof }}%
 {\qed\end{trivlist}}

  \newtheorem{rmk}[theorem]{Remark}
\newenvironment{remark}{\begin{rmk}\em}{\end{rmk}}

  \newtheorem{rmks}[theorem]{Remarks}
\newenvironment{remarks}{\begin{rmks}\em}{\end{rmks}}

  \newtheorem{defn}[theorem]{Definition}
\newenvironment{definition}{\begin{defn}\em}{\end{defn}}

  \newtheorem{eg}[theorem]{Example}
\newenvironment{example}{\begin{eg}\em}{\end{eg}}

  \newtheorem{egs}[theorem]{Examples}
\newenvironment{examples}{\begin{egs}\em}{\end{egs}}


\mathcode`\<="4268  
\mathcode`\>="5269  
\mathcode`\.="313A  
\mathchardef\semicolon="603B 
\mathchardef\gt="313E
\mathchardef\lt="313C

\newcommand{\biimp}
 {\!\Leftrightarrow\!}

\newcommand{\bjg}
 {\mathrel{{\dashv}\,{\vdash}}}

\newcommand{\cod}
 {{\rm cod}}

\newcommand{\dom}
 {{\rm dom}}

\newcommand{\epi}
 {\twoheadrightarrow}

\newcommand{\Ind}[1]
 {{\rm Ind}\hy #1}
 
\newcommand{\imp}
 {\!\Rightarrow\!} 

\newcommand{\mono}
 {\rightarrowtail}

\newcommand{\nml}
 {\triangleleft}

\newcommand{\ob}
 {{\rm ob}}

\newcommand{\op}
 {^{\rm op}}

\newcommand{\pepi}
 {\rightharpoondown\kern-0.9em\rightharpoondown}

\newcommand{\pmap}
 {\rightharpoondown}

\newcommand{\Set}
 {{\bf Set }}

\newcommand{\Sh}
 {{\bf Sh}}

\newcommand{\sh}
 {{\bf sh}}

\newcommand{\Sub}
 {{\rm Sub}}
 
\newcommand{\name}[1]
 {\mbox{$\ulcorner #1 \urcorner$}}

\section{Introduction}

In this paper we give a topos-theoretic interpretation of Priestley duality for distributive lattices, leading to natural analogues of this duality for other categories of partially ordered structures. Specifically, we establish dualities between various categories of ordered structures and categories of Priestley spaces which can be intrinsically characterized through appropriate separation axioms analogously to the case of the classical duality. 

In order to build these `Priestley-type' dualities, we investigate the classical duality from both a topological and an algebraic viewpoint. As it is well-known, topologically the duality is based on the patch topology construction, while algebraically the Boolean algebra of clopen sets of the Priestley space associated to a distributive lattice can be characterized as the free Boolean algebra on it. The unification between the algebraic and the topological formulations of the duality is conveniently provided by the notion of topos; in fact, the toposes involved in Priestley-type dualities admit, on one hand, an algebraic representation (as categories of sheaves on a preordered structure with respect to an appropriate Grothendieck topologies on it) and on the other hand a topological one (as categories of sheaves on natural spectra of the structures, as provided by the techniques of \cite{OC11}).  

Topologically, our `Priestley-type' dualities are built by considering natural spectra for the given partially ordered structures, generating patch-type topologies from them and equipping the resulting spaces with the specialization orderings on the original spectra; algebraically, the dualities are obtained by assigning to any given ordered structure a Boolean algebra which is free on it (in an appropriate sense), equipped with a natural ordering on the points of its spectrum. Specific examples of dualities generated through this method are given in the paper, specifically in section \ref{ex}, and notably include `Priestley-type' dualities for coherent posets, meet-semilattices and disjunctively distributive lattices.  

We also argue more generally that various kinds of free structures provide a natural way for building dualities, to the extent that many free functors admit an inverse defined on an appropriate subcategory. Further illustrations of this phenomenon are provided in the section \ref{ex} of the paper.  

The structure of the paper is as follows. In section \ref{freestmor} we carry our a general analysis of free structures and their construction via syntactic categories, with a particular emphasis on the construction of free Boolean algebras on different kinds of preordered structures through Morleyizations. In section \ref{spatial}, we address the problem of realizing free structures topologically and establish several results which allow us to identify, under natural hypotheses, a free structure on a poset as a structure generated by it inside an appropriate powerset. In section \ref{toposint} we present our topos-theoretic interpretation of Priestley duality, leading to the general method for building `Priestley-type' dualities described in section \ref{general}. The paper ends with a section devoted to concrete examples of dualities generated through our methodologies.  

\section{Free structures and Morleyizations}\label{freestmor}

Let us start with some general remarks about the relationship between free structures and (generalized) syntactic categories.

\subsection{Free structures and syntactic categories}\label{freesyn}

As remarked in \cite{OC11}, the theory of syntactic categories can be profitably applied to the problem of constructing structures presented by generators and relations. In fact, any syntactic category of a given theory $\mathbb T$ can be regarded, in a sense that we shall not make precise in the present paper, as a structure presented by a set of `generators', given by the sorts in the signature of the theory $\mathbb T$, subject to `relations' expressed by the axioms of the theory $\mathbb T$. Conversely, to any structure $\cal C$ one can attach a \emph{canonical signature} $\Sigma_{\cal C}$ to express `relations' holding in the structure, consisting of one sort \name{c} for each element $c$ of $\cal C$ and possibly function or relation symbols whose canonical interpretation in $\cal C$ coincide with specified functions or subsets in $\cal C$ in terms of which the designated `relations' holding in $\cal C$ can be formally expressed; over such a canonical signature one can then write down axioms possibly involving generalized connectives and quantifiers so to obtain a $\cal S$-theory (in the sense of section 8 of \cite{OC10}) $\mathbb T$ whose $\cal S$-syntactic category ${\cal C}^{\cal S}_{\mathbb T}$ can be identified with `the free structure on $\cal C$ subject to the relations $R$', meaning that the $\cal S$-structure $\cal D$ in which the relations $R$ are satisfied naturally correspond to the $\cal S$-homomorphism ${\cal C}^{\cal S}_{\mathbb T} \to {\cal D}$, in a way which can be concretely described as follows. To any $\cal S$-structure $\cal D$ we can canonically associate a $\cal S$-homomorphism ${\cal C}\to {\cal D}$, assigning to any element $c$ of $\cal C$ the interpretation of $\name{c}$ in $\cal D$; in particular we have a canonical $\cal S$-morphism $i:{\cal C} \to {\cal C}^{\cal S}_{\mathbb T}$, in terms of which the universal property of ${\cal C}^{\cal S}_{\mathbb T}$ can be expressed by saying that any $\cal S$-homomorphism $f:{\cal C}\to {\cal D}$ to a $\cal S$-structure $\cal D$ in which the relations $R$ are satisfied can be extended, uniquely up to isomorphism, along the canonical morphism $i$, to a $\cal S$-homomorphism ${\cal C}^{\cal S}_{\mathbb T}\to {\cal D}$.       

As an example of application of this general method, let us consider the problem of constructing, given a small category $\cal C$, a cartesian category $\tilde{{\cal C}}$ with a functor $i:{\cal C}\to \tilde{{\cal C}}$, such that any functor $f:{\cal C}\to {\cal D}$ from $\cal C$ to a cartesian category $\cal D$ can be uniquely extended (up to isomorphism) along $i$ to a cartesian functor $\tilde{f}:\tilde{{\cal C}} \to {\cal D}$. One can write down a cartesian theory $\mathbb T$ over the canonical signature of the category $\cal C$ whose models in any cartesian category can be identified with the functors ${\cal C}\to {\cal D}$ (cf. Example D1.4.8 \cite{El}). Then the category $\tilde{{\cal C}}$, that is the `free cartesian completion' of the category $\cal C$, can be realized as the cartesian syntactic category of the theory $\mathbb T$. Similarly, the `free Boolean completion' of $\cal C$, that is a Boolean coherent category ${\cal B}_{\cal C}$, equipped with a functor $i:{\cal C}\to {\cal B}_{\cal C}$, such that any functor from $\cal C$ to a Boolean coherent category $\cal D$ can be extended along $i$ uniquely to a Boolean coherent functor ${\cal B}_{\cal C} \to {\cal D}$, can be constructed as the first-order syntactic category of the theory $\mathbb T$. Analogously, the regular (resp. coherent, Heyting, geometric, etc.) syntactic category of ${\mathbb T}$ yields the `free regular completion' (resp. `free coherent completion', `free Heyting completion', `free geometric completion', etc.) of the category $\cal C$.

More generally, whenever we can axiomatize, by using a $\cal S$-theory $\mathbb T$ over a canonical signature of a category $\cal C$, a given family $\cal M$ of functors from $\cal C$ to a $\cal S$-category $\cal D$ then the $\cal S$-syntactic category of the theory $\cal S$, together with the canonical functor ${\cal C}\to {\cal C}^{\cal S}_{\mathbb T}$, satisfies the universal property of the `free ${\cal S}$-category on $\cal C$ relative to $\cal M$', that is any functor $f:{\cal C}\to {\cal D}$ from $\cal C$ to a $\cal D$ which belongs to the family $\cal M$ can be extended, uniquely up to isomorphism, to a $\cal S$-functor $\tilde{f}:{\cal C}^{\cal S}_{\mathbb T} \to {\cal D}$.      

Recalling that preorders can be identified with categories in which for any two objects there is at most one arrow from the former to the latter, this methodology can be profitably applied in the context of propositional theories to build ordered algebraic structures presented by generators and relations (cf. \cite{OC11} for a comprehensive treatment of this context); in this paper we shall in particular be concerned with the construction of Boolean algebras presented by generators and relations.   

As much as the logical construction of free structures is important and useful, a problem which frequently arises in practice is that of obtaining explicit descriptions of such structures, descriptions not necessarily of logical nature, but of algebraic, geometric or topological kind, or of whatever other sort. As we already remarked in \cite{OC11}, the methods of Topos Theory can be profitably exploited to obtain such descriptions starting from logical descriptions of the given structures; indeed, any topos which can be naturally attached to the given structure in such a way that the structure can be recovered from it up to isomorphism (for example, a topos from which the structure can be recovered as its full subcategory of $C$-compact objects for a topos-theoretic invariant $C$ - cf. \cite{OC11} for many examples of such situations), admits in general many different representations, which can be obtained by using a variety of topos-theoretic methods; any such representation (of algebraic, resp. geometric or topological) nature will thus produce a corresponding representation (of algebraic, resp. geometric or topological) nature of the given structure. For instance, the classifying topos of a geometric theory can always been built as the topos of sheaves on an appropriate syntactic category of the theory with respect to a natural Grothendieck topology on it, but can also be computed in several different ways (in fact, any Morita-equivalence of the given theory with another one provides a different way of representing the classifying topos, cf. \cite{OC10}), and, at least in the propositional context, these syntactic categories can be recovered from the classifying topos as the full subcategories on the subterminals in the topos which are $C$-compact for a topos-theoretic invariant $C$ (cf. \cite{OC11}). We shall apply these remarks below in connection with the construction of free Boolean algebras on various kinds of partially ordered structures.

\subsection{The Morleyization of a first-order theory}

Let us recall that, given any finitary first-order theory $\mathbb T$ over a signature $\Sigma$, one can define a coherent theory ${\mathbb T}'$ over an extended signature, called the \emph{Morleyization} of $\mathbb T$ (cf. Lemma D1.5.13 \cite{El}), such that for any Boolean coherent category $\cal C$ the category of $\mathbb T$ models in $\cal C$ and elementary morphisms between them is naturally equivalent to the category of models of ${\mathbb T}'$ in $\cal C$. The signature $\Sigma'$ of ${\mathbb T}'$ has, in addition to all the sorts, function symbols and relation symbols of the signature $\Sigma$ of $\mathbb T$, two relation symbols $C_{\phi}\mono A_{1} \cdots A_{n}$ and $D_{\phi}\mono A_{1} \cdots A_{n}$ for each first-order formula $\phi$ over $\Sigma$ (where $A_{1} \cdots A_{n}$ is the string of sorts corresponding to the canonical context of $\phi$), while the axioms of ${\mathbb T}'$ are given by the sequents of the form $C_{\phi} \vdash_{\vec{x}} C_{\psi}$ for any axiom $\phi \vdash_{\vec{x}} \psi$ of $\mathbb T$ together with a set of coherent sequents involving the new relation symbols $C_{\phi}$ and $D_{\phi}$ which ensure that in any model $M$ of ${\mathbb T}'$ in a Boolean coherent category $\cal C$, the interpretation of $C_{\phi}$ coincides with the interpretation of $\phi$ and the interpretation of $D_{\phi}$ coincides with the complement of the interpretation of $\phi$ (cf. p. 859-860 \cite{El} for the details). 

For any finitary first-order theory $\mathbb T$, one can also consider the \emph{classical first-order syntactic category} ${\cal C}^{\textrm{cl-fo}}_{\mathbb T}$ of $\mathbb T$, whose objects are the first-order formulae-in-context over the signature of $\mathbb T$ and whose arrows are the $\mathbb T$-provable (with respect to classical logic) equivalence classes of first-order formulae which are $\mathbb T$-provably functional from the domain to the codomain. This category enjoys the following universal property: for any Boolean coherent category $\cal C$, we have, naturally in $\cal C$, an equivalence of categories
\[
\mathfrak{Bool}({\cal C}^{\textrm{cl-fo}}_{\mathbb T}, {\cal C})\simeq {\mathbb T}\textrm{-mod}_{e}({\cal C}),
\]
where $\mathfrak{Bool}({\cal C}^{\textrm{cl-fo}}_{\mathbb T}, {\cal C})$ denotes the category of Boolean (equivalently, coherent) functors ${\cal C}^{\textrm{cl-fo}}_{\mathbb T} \to {\cal C}$ and natural transformations between them, and ${\mathbb T}\textrm{-mod}_{e}({\cal C})$ denotes the category of models of $\mathbb T$ in $\cal C$ and elementary morphisms between them. This can be easily seen by appropriately modifying the proof of the corresponding result for cartesian theories (as given for example by Theorem D1.4.7 \cite{El}) and observing that the morphisms between the models are elementary since Boolean functors defined on Boolean categories are Heyting functors and hence preserve the interpretations of first-order formulae. 
 
Before proceeding further, we record a couple of easy facts about Morleyizations.

\begin{proposition}\label{morl}
Let $\mathbb T$ be a finitary first-order theory over a signature $\Sigma$ and ${\mathbb T}'$ its Morleyization. Then

\begin{enumerate}[(i)]
\item For any finitary first-order sequent $\phi \vdash_{\vec{x}} \psi$ over $\Sigma$, the sequent is provable in $\mathbb T$ using classical first-order logic if and only if the sequent $C_{\phi} \vdash_{\vec{x}} C_{\psi}$ is provable in ${\mathbb T}'$ using coherent logic; 

\item The classical first-order syntactic category ${\cal C}^{\textrm{cl-fo}}_{\mathbb T}$ of $\mathbb T$ is isomorphic to the coherent syntactic category of ${\mathbb T}'$, and to the classical first-order syntactic category of ${\mathbb T}'$.
\end{enumerate}
\end{proposition}

\begin{proofs}
$(i)$ This follows immediately from the classical completeness theorem for first-order logic and the classical completeness theorem for coherent logic, using the fact that the models of $\mathbb T$ in $\Set$ can be identified with the models of ${\mathbb T}'$ in $\Set$ and that in any such model the interpretation of $\phi$ coincides with the interpretation of $C_{\phi}$.  

$(ii)$ It is immediately verified that every finitary first-order (resp. coherent) formula over the signature of ${\mathbb T}'$ is classically provably equivalent (resp. provably equivalent in coherent logic) to a coherent formula over the signature of ${\mathbb T}'$ (cf. p. 923 \cite{El}); from this it follows at once that the coherent syntactic category of ${\mathbb T}'$ is Boolean and that every morphism between models of ${\mathbb T}'$ in Boolean coherent categories is an elementary morphism. Hence, by the fundamental property of Morleyizations, both the coherent syntactic category of ${\mathbb T}'$ and the classical first-order syntactic category of ${\mathbb T}'$ satisfy the universal property of the category ${\cal C}^{\textrm{cl-fo}}_{\mathbb T}$ with respect to models of $\mathbb T$ in Boolean coherent categories; but this implies, by universality, that these two categories are naturally equivalent to ${\cal C}^{\textrm{cl-fo}}_{\mathbb T}$, as required. 
\end{proofs}

\begin{remarks}
\begin{enumerate}[(a)]
\item A first-order theory is complete in the sense of classical Model Theory (i.e., any first-order sentence over the signature of the theory is either provably false or provably true, but not both) if and only if its Morleyization is complete in the sense of geometric logic (i.e., any geometric sentence over its signature is either provably false or provably true, but not both). Indeed, this immediately follows from the proof of part $(ii)$ of the proposition.

\item From part $(ii)$ of the proposition it immediately follows that two arbitrary first-order theories are equivalent (in the sense that their syntactic Boolean pretoposes are equivalent) if and only if their Morleyizations are Morita-equivalent (i.e. their classifying toposes are equivalent).
\end{enumerate}
\end{remarks}

\subsection{Free Boolean algebras through Morleyizations}\label{freemor}

As we have shown in \cite{OC11}, our general theory of syntactic categories provides a way for building a great variety of preordered structures presented by generators and relations. In this section we shall apply this general technique to construct free Boolean algebras on various kinds of preordered structures, including in particular preorders, meet-semilattices and distributive lattices. Specifically, we shall build such structures as first-order syntactic categories of propositional theories axiomatizing the given class of morphisms from the structure to Boolean algebras, or equivalently as coherent syntactic categories of their Morleyizations (cf. Proposition \ref{morl}). 

The following definition will play a central role in our analysis.

\begin{definition}\label{freedef}
Let $\cal C$ be a structure, $\cal L$ be a category of structures and $\cal M$ be a class of functions from $\cal C$ to structures in $\cal L$. We say that a structure $\cal D$ in $\cal L$, together with a function $i:{\cal C}\to {\cal D}$ in $\cal M$, is the \emph{free $({\cal L}, {\cal M})$-structure} on $\cal C$ if any function $f:{\cal C}\to L$ in $\cal M$ from $\cal C$ to a structure $L$ in $\cal L$ can be uniquely extended via $i$ to an arrow ${\cal D}\to L$ in $\cal L$.
\end{definition} 

It is clearly natural to wonder whether, given $\cal C$, $\cal L$ and $\cal M$ as in the definition, the free $({\cal L}, {\cal M})$-structure on $\cal C$ exists, and if so how it can be built. If the class $\cal L$ can be identified as the class of models of a small ordered algebraic theory $\mathbb A$ and the morphisms in $\cal M$ from $\cal C$ to structures in $\cal L$ can be identified with the models of $\mathbb A$ which satisfy some `relations' written in the canonical signature of $\cal C$ (in the sense of section 8 of \cite{OC11}) then, by Theorem 8.5 \cite{OC11}, the free $({\cal L}, {\cal M})$-structure on $\cal C$ exists and can be built as the $\cal S$-syntactic category of the $\cal S$-theory corresponding to the given set of generators and relations relative to $\mathbb A$ via the method described at pp. 125-126 \cite{OC11}.   

In particular, suppose that $\cal L$ is the category of Boolean algebras. Then $\cal L$ can be seen as the category of models of a small ordered algebraic theory, namely the theory of Boolean algebras, as written in the signature consisting of two constant symbols $0$ and $1$, one unary function symbol $\neg$ and two binary function symbols $\vee$ and $\wedge$. If $\Sigma$ is the canonical signature of $\cal C$ then the $\cal S$-theories over $\Sigma$ corresponding to a set of relations written over $\Sigma$ (via the method of pp. 125-126 \cite{OC11}) are precisely the first-order propositional theories over $\Sigma$. Therefore, the relevant free structures can be built as the syntactic $\cal S$-categories of such theories, that is as the first-order syntactic categories of propositional theories over $\Sigma$, or equivalently as the coherent syntactic categories of their Morleyizations. The use of Morleyizations is important because it allows us to work with coherent theories in place of first-order ones, and hence exploit the theory of classifying toposes to obtain alternative descriptions of the goven syntactic categories, as remarked in section \ref{freesyn} above.    

Let us now give a few examples.

\begin{enumerate}
\item \emph{The free Boolean algebra on a preorder}.\\ Given a preorder $\cal P$, consider the propositional theory $\mathbb P$ over the signature consisting of one relation symbol $R_{a}$ for each element $a$ of $\cal P$ and having as axioms all the sequents of the form
\[
(R_{a} \vdash_{} R_{b})
\]  
for any $a,b\in {\cal P}$ such that $a\leq b$. Clearly, the models of this theory in any Boolean algebra $B$ correspond precisely to the monotone maps ${\cal P}\to B$. The first-order syntactic category ${\cal B}_{\cal P}$ of the theory $\mathbb P$ thus yields a Boolean algebra which is `free on $\cal P$' in the sense that there is a canonical functor $i:{\cal P}\to {\cal B}_{\cal P}$ with the property that any monotone map ${\cal P}\to B$ can be extended uniquely along $i$ to a Boolean algebra homomorphism ${\cal B}_{\cal P}\to B$. 

\item \emph{The free Boolean algebra on a meet-semilattice}.\\ Given a meet-semilattice $\cal M$, consider the theory $\mathbb M$ over the signature consisting of one relation symbol $R_{a}$ for each element $a$ of $\cal M$ and having as axioms all the sequents of the form
\[
(\top \vdash R_{1_{{\cal M}}}),
\]
\[
(R_{a} \wedge R_{b} \dashv\vdash_{} R_{a\wedge b})
\]  
for any $a,b\in {\cal M}$, where $1_{\cal M}$ denotes the top element of $\cal M$ and $\wedge$ denotes the meet operation in $\cal M$. 

The models of $\mathbb M$ in a Boolean algebra $B$ clearly correspond to the meet-semilattice homomorphisms ${\cal M}\to B$; therefore, the first-order syntactic category of $\mathbb M$ provides the free Boolean algebra on $\cal M$; that is, we have a canonical functor $i:{\cal M}\to {\cal B}_{\cal M}$ with the property that any meet-semilattice homomorphism ${\cal M}\to B$ to a Boolean algebra $B$ can be extended uniquely along $i$ to a Boolean algebra homomorphism ${\cal B}_{\cal M}\to B$.   

\item \emph{The free Boolean algebra on a distributive lattice}.\\ Given a distributive lattice $\cal D$, consider the theory $\mathbb D$ over the signature having one $0$-ary relation symbol $R_{d}$ for each element $d\in D$, and the following axioms:
\[
(\top \vdash R_{1_{\cal D}});
\]
\[
(R_{0_{\cal D}} \vdash \bot);
\]
\[
(R_{a\wedge b} \dashv\vdash R_{a} \wedge R_{b})
\] 
for any $a, b$ in $\cal D$; 
\[
(R_{a \vee b} \dashv\vdash R_{a} \vee R_{b})
\]
for any $a, b \in {\cal D}$ (where $1_{\cal D}$ denotes the top element of $\cal D$, $0_{\cal D}$ denotes the bottom element of $\cal D$, $\wedge$ denotes the meet operation in $\cal D$ and $\vee$ denotes the join operation in $\cal D$).

For any Boolean algebra $B$, the models of ${\cal D}$ in $B$ correspond precisely to the distributive lattice homomorphisms ${\cal D} \to B$. The first-order syntactic category of $\mathbb D$ thus provides the free Boolean algebra on $\cal D$, in the sense that we have a canonical functor $i:{\cal D}\to {\cal B}_{\cal D}$ with the property that any distributive lattice homomorphism ${\cal D}\to B$ to a Boolean algebra $B$ can be extended uniquely along $i$ to a Boolean algebra homomorphism ${\cal B}_{\cal D}\to B$.     

\end{enumerate}

\subsubsection{The free Boolean algebra on a distributive lattice}

In this section we give, as an example of our technique of using classifying toposes for obtaining concrete descriptions of structures presented by generators and relations, an explicit description of the free Boolean algebra on a distributive lattice.

By the discussion in the last section, the free Boolean algebra on a distributive lattice $\cal D$ can be identified with the classical first-order syntactic category of the theory $\mathbb D$ introduced in section \ref{freemor} above; this category can in turn be identified, by Proposition \ref{morl}, with the coherent syntactic category of the Morleyization of $\mathbb D$. In fact, the Morleyization of $\mathbb D$ admits in this case a simpler description, that is we can alternatively use, in place of it, a simpler theory having an equivalent coherent syntactic category. This theory, which we call ${\cal S}_{\cal D}$, can be described as follows. Its signature has, in addition to all the $0$-ary relation symbols of the theory $\mathbb D$, one $0$-ary relation symbol $R^{\ast}$ for each ($0$-ary) relation symbol $R$ over the signature of $\mathbb D$, and, in addition to the axioms of $\mathbb D$, the following axioms:
\[ 
(\top \vdash R_{d} \vee R_{d}^{\ast})
\]      
for any $d\in D$, and 
\[
(R_{d} \wedge R_{d}^{\ast} \vdash \bot)
\]
for any $d\in D$.

The coherent syntactic category ${\cal C}_{{\cal S}_{\cal D}}$ of this theory, together with the obvious map ${\cal D}\to {\cal C}_{{\cal S}_{\cal D}}$ (sending any element $d\in {\cal D}$ to the equivalence class of the formula $R_{d}$) thus yields the free Boolean algebra on the distributive lattice $\cal D$. This provides a logical description of this structure which, in particular, establishes its existence independently of any non-constructive assumptions (such as for example weak forms of the axiom of choice necessary to establish spatial representations for the structures in question); in fact, it is well-known that the Stone representation of any distributive lattice as a ring of sets gives naturally rise to a related description of the free Boolean algebra on it (cf. also below), but such spatial representations are not guaranteed to exist in a fully constructive framework. Concerning this, it should be mentioned that \cite{CL} gives an elegant logical description of the free Boolean algebra on a distributive lattice in terms of entailment relations, while \cite{peremans} provides a more algebraic, but still intrinsically logical in character, description of the same structure. In this section we present an alternative algebraic description of this structure, directly involving the elements and operations on the lattice and not relying on `deduction rules' of any sort. Similarly to how we `computed' the free frame on a complete join-semilattice, we achieve such a description by investigating the classifying topos of the theory ${\cal C}_{{\cal S}_{\cal D}}$ from a `semantical' point of view. 

As shown in \cite{OC6}, a convenient way for obtaining a `semantic' representation for the classifying topos of a geometric theory consists in regarding the theory as a quotient of a theory of presheaf type; this produces a representation of the classifying topos of the theory as a topos of sheaves on the opposite of the category of finitely presentable models of the relevant theory of presheaf types with respect to a Grothendieck topology which can be directly described in terms of the axioms of the theory (cf. \cite{OC6}) for the details of the general theory). In our case, a convenient choice is that of regarding the theory ${\cal C}_{{\cal S}_{\cal D}}$ as a quotient of the empty theory over its signature (note however that this is by no means the only possibility, and that different choices would produce different representations for the same classifying topos). This leads to a description of the free Boolean algebra on $D$ much in the same spirit as that of the free frame on a complete-join semilattice obtained in \cite{OC11}. 

Concerning notation, we denote the disjoint union of $\cal D$ with itself by ${\cal D}_{db}$; for any element $d\in D$, considered as an element of the `first copy of $D$' in ${\cal D}_{db}$, we denote by $d^{\ast}$ the corresponding element in the second copy of $\cal D$ in ${\cal D}_{db}$. 

\begin{theorem}\label{freeBoolean}
The free Boolean algebra on a distributive lattice $\cal D$ can be realized as the set $B({\cal D})$ of compact elements of the frame of upward closed subsets ${\cal I}$ of the set $\mathscr{P}_{fin}({\cal D}_{db})$ of finite subsets of ${\cal D}_{db}$ (with respect to the subset-inclusion ordering on $\mathscr{P}_{fin}({\cal D}_{db})$) satisfying the following properties (endowed with the subset-inclusion ordering):

\begin{enumerate}

\item for any subset $U\in \mathscr{P}_{fin}({\cal D}_{db})$, $U\cup \{1_{\cal D}\}\in {\cal I}$ implies $U\in {\cal I}$;

\item for any subset $U\in \mathscr{P}_{fin}({\cal D}_{db})$, $0_{\cal D}\in U$ implies $U\in {\cal I}$;

\item for any elements $a,b\in {\cal D}$, if $U\cup \{a, a \vee b\}\in {\cal I}$ and $U\cup \{b, a \vee b\}\in {\cal I}$ then $U\cup \{a \vee b \}\in {\cal I}$;

\item for any elements $a,b\in {\cal D}$, if $U\cup \{a,b, a\wedge b\}\in {\cal I}$ then $U\cup \{a, b\}\in {\cal I}$;

\item for any $d\in {\cal D}$ and any subset $U\in \mathscr{P}_{fin}({\cal D}_{db})$, $U\cup \{d, d^{\ast}\}\in {\cal I}$;

\item for any $d\in {\cal D}$ and any subset $U\in \mathscr{P}_{fin}({\cal D}_{db})$, if $U\cup \{d\}\in {\cal I}$ and $U\cup \{d^{\ast}\}\in {\cal I}$ then $U\in {\cal I}$, 
\end{enumerate}

endowed with the induced ordering, with the universal map ${\cal D}\to B({\cal D})$ being the function sending any element $d\in {\cal D}$ to the smallest element of $B({\cal D})$ containing the set $\{d\}$ as one of its elements.   
\end{theorem}\qed    

\begin{proofs}
This follows immediately from the explicit description of the classifying topos $\Sh(\mathscr{P}_{fin}({\cal D}_{db})^{\textrm{op}}, J)$ of the theory ${\cal C}_{{\cal S}_{\cal D}}$ in terms of its axioms, and from the remark that the coherent syntactic category of ${\cal C}_{{\cal S}_{\cal D}}$ can be recovered as the set of compact elements of the frame of subterminals of the classifying topos of ${\cal C}_{{\cal S}_{\cal D}}$ (that is, of the frame of $J$-ideals on $\mathscr{P}_{fin}({\cal D}_{db})^{\textrm{op}}$, cf. \cite{OC11}). 
\end{proofs}

\subsection{Spatial realization of free structures}\label{spatial}

If a given partially ordered structure can be represented as a substructure of a powerset then it is natural to wonder if structures which are free on it (in the sense of Definition \ref{freedef} above) can also be naturally realized as substructures of that powerset.

Suppose that $\cal L$ is a category of structures and homomorphisms between them and that every powerset is equipped with (finitary) function and relation symbols which make it into a structure in $\cal L$. Given a set $X$ and a subset ${\cal C}\subseteq \mathscr{P}(X)$, there always exists a smallest substructure of $\mathscr{P}(X)$ in $\cal L$ containing $\cal C$, called the substructure of $\mathscr{P}(X)$ generated by $\cal C$ (cf. pp. 7-8 \cite{Hodges}); we shall denote it by $G^{X}_{({\cal C}, {\cal L})}$.

The following proposition exhibits a relationship between the concept of substructure and that of free structure.

\begin{proposition}\label{proplemma}
Let $\cal C$ be a structure, $\cal L$ be a category of structures and $\cal M$ be a class of functions from $\cal C$ to structures in $\cal L$. Suppose that ${\cal C}\subseteq {\cal D} \subseteq \mathscr{P}(X)$ and that the inclusion $m:{\cal C}\hookrightarrow {\cal D}$ identifies $\cal D$ as the free $({\cal L}, {\cal M})$-structure on $\cal C$. Then, if the inclusion ${\cal C}\hookrightarrow G^{X}_{({\cal C}, {\cal L})}$ is a morphism in $\cal M$, we have ${\cal D}=G^{X}_{({\cal C}, {\cal L})}$, that is $G^{X}_{({\cal C}, {\cal L})}$ is the free $({\cal L}, {\cal M})$-structure on $\cal C$.  
\end{proposition}

\begin{proofs}
By the universal property of $G^{X}_{({\cal C}, {\cal L})}$, we have an inclusion $j:G^{X}_{({\cal C}, {\cal L})}\subseteq {\cal D}$, so it remains to prove the converse inclusion. Since the inclusion $i:{\cal C}\hookrightarrow G^{X}_{({\cal C}, {\cal L})}$ is a morphism in $\cal M$, there exists an arrow $r:{\cal D}\to G^{X}_{({\cal C}, {\cal L})}$ in $\cal L$ such that $r\circ m=i$. In order to prove the inclusion ${\cal D} \hookrightarrow G^{X}_{({\cal C}, {\cal L})}$, we will show that $t\circ r=u$, where $t$ is the inclusion $G^{X}_{({\cal C}, {\cal L})} \hookrightarrow \mathscr{P}(X)$ and $u$ is the inclusion ${\cal D} \hookrightarrow \mathscr{P}(X)$. Notice that $u\circ j=t$, since all these functions are inclusions. Now, by the uniqueness of the arrow $k:{\cal D}\to {\cal D}$ such that $k\circ m=m$ given by the universal property of the free $({\cal L}, {\cal M})$-structure on $\cal C$, we have that $j\circ r=1_{\cal D}$, whence $u=u\circ 1_{\cal D}=u\circ j\circ r=t\circ r$, as required. 
\end{proofs}

\begin{remark}
In general, $G^{X}_{({\cal C}, {\cal L})}$ is \emph{not} the free $({\cal L}, {\cal M})$-structure on $\cal C$. Take for example ${\cal C}=\mathscr{P}(X)$; $\cal C$, regarded as a meet-semilattice, coincides with the Boolean algebra generated by itself, but it is not the free Boolean algebra on it (as a meet-semilattice), since Boolean algebra homomorphisms of Boolean algebras do not in general coincide with meet-semilattice homomorphisms of their underlying meet-semilattices. Anyway, as we shall see below, for a great variety of structures $\cal C$ there is a natural choice of a set $X$ such that $\cal C$ can be embedded as a substructure of the powerset $\mathscr{P}(X)$ and the free $({\cal L}, {\cal M})$-structure on $\cal C$ can be identified with the substructure $G^{X}_{({\cal C}, {\cal L})}$ in $\cal L$ of $\mathscr{P}(X)$ generated by $\cal C$.
\end{remark}

The following result provides a set of natural contexts in which a free structure on a given one can be identified as a substructure generated by it.

\begin{theorem}\label{prgen}
Let $\cal C$ be a preordered structure, $\cal L$ be a category of preordered structures and $\cal M$ be a class of functions from $\cal C$ to structures in $\cal L$. Let assume that all the powersets can be made into structures in $\cal L$ and that for any function $f:X\to Y$ the induced function $\mathscr{P}(f)=f^{-1}:\mathscr{P}(Y)\to \mathscr{P}(X)$ is an arrow in $\cal L$. Let us moreover assume that, given a function from $\cal C$ to a structure in $\cal L$, if the composition of it with each of the arrows belonging to a jointly injective family of arrows in $\cal L$ belongs to $\cal M$ then the function itself belongs to $\cal M$. 
Let the map $m:{\cal C}\to L$ identify $L$ as the free $({\cal L}, {\cal M})$-structure on $\cal C$. Suppose that there exists a subcanonical Grothendieck topology $J_{\cal C}$ on $\cal C$ such that the morphisms ${\cal C}\to \{0,1\}$ in $\cal M$ coincide precisely with a designated set $X_{\cal C}$ of jointly conservative flat $J_{\cal C}$-continuous functors ${\cal C}\to \{0,1\}$, and that there is a subcanonical Grothendieck topology $K_{L}$ on $L$ such that the morphisms $L\to \{0,1\}$ in $\cal L$ which extend the morphisms ${\cal C}\to \{0,1\}$ in $\cal M$ via $m$ coincide with a designated set $Y_{L}$ of jointly conservative flat $K_{L}$-continuous functors $L \to \{0,1\}$. Then $L$ can be identified with the $\cal L$-substructure of $\mathscr{P}(X_{\cal C})$ generated by $\cal C$, regarded as a subset of $\mathscr{P}(X_{\cal C})$ via the composite map ${\cal C}\mono Id_{J_{\cal C}}({\cal C})\cong {\cal O}(X_{\cal C})\subseteq \mathscr{P}(X_{\cal C})$ (where ${\cal C}\mono Id_{J_{\cal C}}({\cal C})$ is the canonical embedding of $\cal C$ into the frame $Id_{J_{\cal C}}({\cal C})$ of $J_{\cal C}({\cal C})$-ideals on $\cal C$ and $X_{\cal C}$ is endowed with the topology induced by that of the space of points of the topos $\Sh({\cal C}, J_{\cal C})$). We have a geometric morphism
\[
\Sh(L, K_{L})\simeq \Sh(Y_{L})\to \Sh(X_{\cal C})\simeq \Sh({\cal C}, J_{\cal C}),
\]
where $Y_{L}$ is endowed with the topology induced by that of the space of points of the topos $\Sh(L, K_{L})$.

If moreover the $\cal L$-substructure of $\mathscr{P}(X_{\cal C})$ generated by $\cal C$ is contained in ${\cal O}(X_{\cal C})$ then we have an equivalence of toposes
\[
\Sh(L, K_{L})\simeq \Sh(Y_{L})\simeq \Sh(X_{\cal C}) \simeq \Sh({\cal C}, J_{\cal C}).
\]
\end{theorem}

\begin{proofs}
Since $K_{L}$ is subcanonical and the topos $\Sh(L, K_{L})$ has enough points, we have an embedding $k:L\to \mathscr{P}(Y_{L})$. Similarly, we have an embedding $h:{\cal C}\to \mathscr{P}(X_{\cal C})$.
Since, by our hypotheses, the arrows ${\cal C}\to \{0,1\}$ in $\cal M$ correspond bijectively (by composition with $m$) to the arrows $L\to \{0,1\}$ in $Y_{L}$, we have a bijection $u:Y_{L}\to X_{\cal C}$ such that the diagram
\[  
\xymatrix {
{\cal C} \ar[d]^{m} \ar[rr]^{h} & & \mathscr{P}(X_{\cal C}) \ar[d]^{\mathscr{P}(u)}   \\
L \ar[rr]^{k} & & {\mathscr{P}}(Y_{L})}
\]
commutes. Notice that the commutativity of this diagram implies in particular that the arrow $m$ is injective and that therefore the free $({\cal L}, {\cal M})$-structure on $\cal C$ can be realized (up to isomorphism) as a subset of $\mathscr{P}(X_{\cal C})$ containing $\cal C$. In order to conclude our thesis by using Proposition \ref{proplemma}, it only remains to verify that the inclusion ${\cal C} \hookrightarrow G^{X_{\cal C}}_{({\cal C}, {\cal L})}$ is a morphism in $\cal M$. To this end, we observe the following facts. The arrows of the form $\mathscr{P}(\overline{x}):\mathscr{P}(X_{\cal C}) \to \{0,1\}\cong \mathscr{P}(\{\ast\})$, where $\overline{x}:\{\ast\}\to X_{\cal C}$ is the function defined by $\overline{x}(\ast)=x$, belong to $\cal L$ (by our hypotheses) and are jointly injective (for $x\in X_{\cal C}$); therefore, since for every $x\in X_{\cal C}$ $\mathscr{P}(\overline{x})\circ h$ belongs to $\cal M$, the function $h:{\cal C}\to \mathscr{P}(X_{\cal C})$ belongs to $\cal M$. From this it follows in turn, by invoking our hypotheses again, that the factorization ${\cal C} \hookrightarrow G^{X_{\cal C}}_{({\cal C}, {\cal L})}$ of $h$ across the inclusion $G^{X_{\cal C}}_{({\cal C}, {\cal L})} \hookrightarrow \mathscr{P}(X_{\cal C})$ (which is an arrow of $\cal L$ by definition of $\cal L$-substructure of $\mathscr{P}(X_{\cal C})$ generated by $\cal C$) belongs to $\cal M$, as required.   

The commutativity of the square above implies that the map $u:Y_{L}\to X_{\cal C}$ is continuous if $Y_{L}$ (resp. $X_{\cal C}$) is endowed with the topology induced by that of the space of points of the topos $\Sh(L, K_{L})$ (resp. of the topos $\Sh({\cal C}, J_{\cal C})$), whence it induces a geometric morphism
\[
\tilde{u}:\Sh(L, K_{L})\simeq \Sh(Y_{L})\to \Sh(X_{\cal C})\simeq \Sh({\cal C}, J_{\cal C}).
\] 
The last part of the theorem can be proved as follows. If $L$ corresponds, under the bijection $\mathscr{P}(u)$, to a subset of $\mathscr{P}(X_{\cal C})$ contained in ${\cal O}(X_{\cal C})$ then the frame of open sets ${\cal O}(Y_{L})$ of $Y_{L}$ corresponds under $\mathscr{P}(u)$ to the set ${\cal O}(X_{\cal C})$ of open sets of the space $X_{\cal C}$ (since $L$ constitutes a basis for $Y_{L}$). Hence the spaces $X_{\cal C}$ and  $Y_{L}$ are homeomorphic (under the bijection $u$) whence the toposes $\Sh(X_{\cal C})$ and $\Sh(Y_{L})$ are equivalent (under the geometric morphism $\tilde{u}$ defined above). 
\end{proofs}

\begin{remarks}\label{remvar}
\begin{enumerate}[(a)]
\item Under the hypotheses stated in the first paragraph of the theorem, a sufficient set of conditions for the remaining hypotheses to be satisfied is the following: there is a subcanonical Grothendieck topology $J_{\cal C}$ on $\cal C$ and for any $L\in {\cal L}$ a subcanonical Grothendieck topology $K_{L}$ on $L$ such that:

\begin{enumerate}[(i)]
\item For any $L, L'\in {\cal L}$, the morphisms $L\to L'$ in $\cal L$ coincide precisely with the morphisms of sites $(L, K_{L})\to (L', K_{L'})$;

\item For any $L\in {\cal C}$, the morphisms ${\cal C} \to L$ in $\cal M$ coincide precisely with the morphisms of sites $({\cal C}, J_{\cal C})\to (L, K_{L})$;

\item The toposes $\Sh({\cal C}, J_{\cal C})$ and $\Sh(L, K_{L})$ (for $L$ in $\cal L$) all have enough points. 
\end{enumerate} 
Indeed, under conditions $(ii)$ and $(iii)$, the morphisms ${\cal C}\to \{0,1\}$ in $\cal M$ coincide exactly with the morphisms of sites $({\cal C}, J_{\cal C})\to (\{0,1\}, K_{\{0,1\}})$ and hence with the $J_{\cal C}$-continuous flat functors ${\cal C}\to \{0,1\}$ (notice that, it being subcanonical, $K_{\{0,1\}}$ necessarily coincides with the Grothendieck topology $K$ on $\{0,1\}$ such that $K(0)$ consists of the empty sieve and the maximal one, and $K(1)$ consists of the maximal sieve only); and these flat functors are jointly conservative if and only if the topos $\Sh({\cal C}, J_{\cal C})$ has enough points. Similarly, it follows from $(i)$ and $(iii)$ that, for any $L\in {\cal L}$, the morphisms $L\to \{0,1\}$ in $\cal L$ correspond precisely to the points of the topos $\Sh(L, K_{L})$ and are jointly conservative for it.

\item We shall apply the theorem only in the particular case of the assumptions of part $(a)$ above where the class $\cal M$ can be identified with the class of geometric morphisms $\Sh(L, K_{L})\to \Sh({\cal C}, J_{\cal C})$ (that is, with the class of models of the theory of $J_{\cal C}$-prime filters on $\cal C$ in the topos $\Sh(L, K_{L})$, cf. \cite{OC11}). Notice in passing that under these hypotheses both the spaces $X_{\cal C}$ and $Y_{L}$ are sober, since they are homeomorphic to spaces of points of localic toposes.    

\item The morphism $m:{\cal C}\to L$ realizing $L$ as the free $({\cal L}, {\cal M})$-structure on $\cal C$ is a morphism of sites $({\cal C}, J_{\cal C})\to (L, K_{L})$. Indeed, the geometric morphism $\tilde{u}:\Sh(Y_{L})\to \Sh(X_{{\cal C}})$ defined in the proof of the theorem satisfies the property that its inverse image $\tilde{u}^{\ast}$ restricts to the map $m:{\cal C}\to L$, and hence $m$ must be a morphism of sites $({\cal C}, J_{\cal C})\to (L, K_{L})$ inducing $\tilde{u}$ (cf. Lemma C2.3.8 \cite{El}). 
\end{enumerate}
\end{remarks}

\section{A topos-theoretic look at Priestley duality}\label{toposint}

In this section we provide a topos-theoretic interpretation of Priestley duality for distributive lattices; this interpretation will pave the way for the general setup for building `Priestley-type dualities' described in section \ref{general} below.

Before reviewing the classical duality, we shall embark in a general analysis of the concept of preordered topological space, which plays a crucial role in the duality.

\subsection{Preordered topological spaces}\label{pts}

We define a \emph{preordered topological space} as a triple $(X, \tau, \leq)$, where $X$ is a set, $\tau$ is a topology on $X$ and $\leq$ is a preorder relation on $X$. Preordered topological spaces form a category, which we denote $\textbf{PTop}$, whose arrows $(X, \tau, \leq)\to (X', \tau', \leq')$ are the maps $f:X\to Y$ which are preorder-preserving (i.e., such that for any $x,y\in X$, $x\leq y$ implies $f(x)\leq' f(y)$) and continuous (i.e., such that for any open set $V\in \tau'$, $f^{-1}(V)\in \tau$); composition and identities in $\textbf{PTop}$ are defined by composing the underlying functions set-theoretically.

We have a functor $i_{P}:\textbf{Top}\to \textbf{PTop}$, sending a topological space $(X, \tau)$ to the triple $(X, \tau, \leq_{X})$, where $\leq_{X}$ is the \emph{specialization preorder} on $X$, that is the preorder relation on $X$ defined by: for any $x,y\in X$, $x\leq_{X} y$ if and only if for every $U\in \tau$, $x\in U$ implies $y\in U$. Indeed, any continuous map $f:(X, \tau)\to (X', \tau')$ is preorder-preserving with respect to the specialization preorders on $X$ and $X'$. We can define a functor $r_{P}:\textbf{PTop} \to \textbf{Top}$ which is left adjoint to $i_{P}$, as follows.
For any object $(X, \tau, \leq)$ of $\textbf{PTop}$, we set $r_{P}((X, \tau, \leq))$ equal to the pair $(X, \tau_{\leq})$, where $\tau_{\leq}$ is the topology on $X$ whose open sets are exactly the open sets of $\tau$ which are $\leq$-upper sets, and for any arrow $f:(X, \tau, \leq) \to (X', \tau', \leq')$ in $\textbf{PTop}$, we set $r_{P}(f)$ equal to $f:(X, \tau_{\leq}) \to (X', \tau'_{\leq'})$. Notice that this is indeed a continuous map of topological spaces since the fact that $f$ is order-preserving ensures that for every open set $U'$ of $\tau'$, if $U'$ is a $\leq'$-upper set then $f^{-1}(U)$ is a $\leq$-upper set.

Note that this definition also makes sense more generally for an arbitrary binary relations in place of preorders, but in fact there is no loss of generality in supposing $\leq$ to be preorder since for any binary relation symbol $R$ on $X$, $r_{P}((X, \tau, R))=r_{P}((X, \tau, \dot{R}))$, where $\dot{R}$ is the reflexive and transitive closure of $R$.

\begin{proposition}\label{proptop}
With the notation above, the functor $i_{P}:\textbf{Top}\to \textbf{PTop}$ is right adjoint to the functor $r_{P}:\textbf{PTop} \to \textbf{Top}$, and identifies $\textbf{Top}$ with a full reflective subcategory of $\textbf{PTop}$.   
\end{proposition}  
 
\begin{proofs}
We need to show that we have a bijective correspondence between the arrows $r_{P}((X, \tau, \leq)) \to (Y, \tau')$ in $\textbf{Top}$ and the arrows $(X, \tau, \leq) \to (Y, \tau', \leq_{Y})$ in $\textbf{PTop}$, naturally in any object $(X, \tau, \leq)$ of $\textbf{PTop}$ and any object $(Y, \tau')$ of $\textbf{Top}$. Given an arrow $f:r_{P}((X, \tau, \leq))=(X, \tau_{\leq}) \to (Y, \tau')$ in $\textbf{Top}$, we associate to it the arrow $\tilde{f}:(X, \tau, \leq) \to (Y, \tau', \leq_{Y})$ of $\textbf{PTop}$ whose underlying function is $f$. Conversely, given an arrow $g:(X, \tau, \leq) \to (Y, \tau', \leq_{Y})$ in $\textbf{PTop}$, we associate to it the arrow $\hat{g}:r_{P}((X, \tau, \leq))=(X, \tau_{\leq}) \to (Y, \tau')$ in $\textbf{Top}$ whose underlying function is $g$. 

In order to show that these assignments are well-defined, it suffices to check that the following conditions on a function $f:X\to Y$ are equivalent:

\begin{enumerate}[(i)]
\item for any $x,x'\in X$, if $x\leq x'$ then $f(x)\leq_{Y} f(x')$;

\item for any open set $V\in \tau'$, $f^{-1}(V)$ is a $\leq$-upper set in $X$. 
\end{enumerate}  
To show that $(i)\imp (ii)$, we suppose that $V$ is an open set in $\tau'$ and prove that, for any $x, x'\in X$ such that $x\leq x'$, if $x\in f^{-1}(V)$ then $x'\in f^{-1}(V)$. If $x\leq x'$ then $f(x)\leq_{Y} f(x')$ and therefore $f(x)\in V$ implies $f(x')\in V$ by definition of specialization preorder $\leq_{Y}$, as required. 

To show that $(ii)\imp (i)$, suppose that $x\leq x'$ in $X$. To prove that $f(x)\leq_{Y} f(x')$ we have to verify that for every $V\in \tau'$, $f(x)\in V$ (equivalently, $x\in f^{-1}(V)$) implies that $f(x')\in V$ (equivalently, $x'\in f^{-1}(V)$); but this follows immediately from the fact that $f^{-1}(V)$ is a $\leq$-upper set.

Clearly, the assignments $f\to \tilde{f}$ and $g\to \hat{g}$ are inverse to each other and natural in $(X, \tau, \leq)\in \textbf{PTop}$ and $(Y, \tau')\in \textbf{Top}$.

Moreover, the facts that $i_{P}$ is full and that $r_{P} \circ i_{P}$ is isomorphic to the identity functor on $\textbf{Top}$ are clear. The proof of the proposition is therefore complete.
\end{proofs} 

Let $\textbf{Pro}$ be the category of preorders and monotone maps between them. We have a functor $L_{\textbf{Pro}}:\textbf{Pro} \to \textbf{PTop}$ sending a preorder $(P, \leq)$ to the triple $(P, \tau_{P}, \leq)$, where $\tau_{P}$ is the Alexandrov topology on $P$, and sending a monotone map $f:(P, \leq) \to (Q, \leq')$ to the arrow $f:(P, \tau_{P}, \leq) \to (Q, \tau_{Q}, \leq')$ of $\textbf{PTop}$. We can define a functor $R_{\textbf{Pro}}:\textbf{PTop}\to \textbf{Pro}$ which is right adjoint to $L_{\textbf{Pro}}$, as follows. $R_{\textbf{Pro}}$ sends an object $(X, \tau, \leq)$ of $\textbf{PTop}$ to $(X, \dot{\leq})$, where $\dot{\leq}$ is the preorder on $X$ given by the intersection between $\leq$ and the specialization preorder on $X$ induced by the topology $\tau$, and an arrow $f:(X, \tau, \leq)\to (Y, \tau', \leq')$ to the monotone map $f:(X, \dot{\leq})\to (Y, \dot{\leq'})$. 

\begin{proposition}
With the notation above, the functor $L_{\textbf{Pro}}:\textbf{Pro} \to \textbf{PTop}$ is left adjoint to the functor $R_{\textbf{Pro}}:\textbf{PTop} \to \textbf{Top}$, and identifies $\textbf{Pro}$ with a full coreflective subcategory of $\textbf{PTop}$.   
\end{proposition}  
 
\begin{proofs}
We have to prove that there is a bijective correspondence between the arrows 
\[
L_{\textbf{Pro}}(P, \leq)=(P, \tau_{P}, \leq) \to (X, \tau, \leq)
\]
in $\textbf{PTop}$ and the arrows 
\[
(P, \leq) \to (X, \dot{\leq})
\]
in $\textbf{Pro}$, naturally in any object $(P, \leq)$ of $\textbf{Pro}$ and any object $(X, \tau, \leq)$ of $\textbf{PTop}$. 

Given an arrow $f:(P, \tau_{P}, \leq) \to (X, \tau, \leq)$ in $\textbf{PTop}$, we associate to it the arrow $\tilde{f}:(P, \leq)\to (X, \dot{\leq})$ whose underlying function is $f$; conversely, given an arrow $g:(P, \leq)\to (X, \dot{\leq})$ in $\textbf{Pro}$, we associate to it the arrow $\hat{g}:(P, \tau_{P}, \leq) \to (X, \tau, \leq)$ in $\textbf{PTop}$ whose underlying function is $g$. In order to show that these assignments are well-defined, it suffices to verify that for any monotone map $f:(P, \leq)\to (X, \leq)$, $f$ is a monotone map $(P, \leq)\to (X, \dot{\leq})$ if and only if it is a continuous map $(P, \tau_{P})\to (X, \tau)$. This follows by similar arguments as those in the proof of Proposition \ref{proptop}.   

It is also immediate to see that the functor $L_{\textbf{Pro}}$ is full and that the composite functor  $R_{\textbf{Pro}} \circ L_{\textbf{Pro}}$ is isomorphic to the identity functor on $\textbf{Pro}$. This completes the proof of the proposition. 
\end{proofs}

By composing the adjunctions obtained in the propositions above, we recover the well-known adjunction between the category of preorders and the category of topological spaces. We can represents the results obtained in this section by means of the following diagram.

\[  
\xymatrix {
&  \textbf{PTop} \ar@<0.9ex>[dl]^{r_{P}} \ar@<1.4ex>[dr]^{R_{\textbf{Pro}}} &   \\
\textbf{Top} \ar@<1.4ex>[ur]^{i_{P}}_{\sevdash} \ar@<0.6ex>[rr]_{\downvdash} & & \textbf{Pro} \ar@<0.9ex>[ll] \ar@<0.9ex>[ul]^{L_{\textbf{Pro}}}_{\swvdash}}
\]

\subsection{Review of Priestley duality}

Priestley duality for distributive lattices is a categorical duality between the category of distributive lattices and the category of \emph{Priestley spaces}. Via this duality, a distributive lattice $D$ corresponds to the ordered topological space $P_{D}$ obtained by equipping the set ${\cal F}_{D}$ of prime filters on $D$ with the \emph{patch topology} (i.e., the topology having as a sub-basis the collection of the sets of the form $\{P\in {\cal F}_{D} \textrm{ | } d\in P\}$ for $d\in D$ and their complements) and with the specialization order $\leq$ on ${\cal F}_{D}$ induced by the coherent topology on ${\cal F}_{D}$ (i.e., the topology having as a basis the collection of sets of the form $\{P\in {\cal F}_{D} \textrm{ | } d\in P\}$). The assignment $D\to P_{D}$ can be made functorial as follows: any morphism $D\to D'$ of distributive lattices induces an order-preserving continuous map $P_{D'}\to P_{D}$; therefore, if we denote by $\textbf{PTop}$ the category of ordered topological spaces we have a functor $P:\textbf{DLat}\to \textbf{PTop}$, to which we shall refer as the \emph{Priestley functor}.  

Any distributive lattice $D$ can be recovered from the associated Priestley space as the poset of clopen $\leq$-upper sets in it; in fact, this assignment defines a functor from the category of Priestley spaces to the opposite of the category of distributive lattices which yields the other half of Priestley duality. 

The ordered topological spaces which are, up to isomorphism, in the image of the Priestley functor are called Priestley spaces; notably, these spaces admit natural intrinsic topological characterizations (for instance, as the compact ordered topological spaces satisfying the Priestley separation axiom or, alternatively, as the ordered topological spaces having a sub-basis of clopen up-sets and clopen down-sets).
 
Priestley duality admits a purely topological interpretation as a categorical equivalence between the category of coherent spaces and the category of Priestley spaces (cf. \cite{cornish}); in fact, this equivalence can be obtained by composing Stone duality between coherent spaces and distributive lattices and Priestley duality between distributive lattices and Priestley spaces. Specifically, a coherent space $(X, \tau)$ corresponds to the Priestley space $(X, \tau', \leq)$, where $\tau'$ is the topology on $X$ having as a sub-basis the set of compact open sets of $(X, \tau)$ and their complements, while a Priestley space $U$ corresponds to the coherent space on the set $X$ whose open sets are exactly the open upper sets of $U$.

We can also naturally look at Priestley duality from an algebraic viewpoint. As remarked in \cite{stone} (cf. the Exercise at p. 73), the algebra of clopen subsets of the Priestley space associated to a distributive lattice $D$ via Priestley duality can be characterized as the free Boolean algebra on $D$. In fact, this characterization can be obtained as an immediate consequence of our Theorem \ref{prgen} in light of the fact that the free Boolean algebra on any distributive lattice exists (for example, by our syntactic method for constructing structures presented by generators and relations, cf. section 8 of \cite{OC11}).    

In the following section, we introduce a general abstract framework for interpreting Priestley duality, leading to `Priestley-type' dualities for various classes of partially ordered structures other than distributive lattices. Our framework can be equivalently presented in the language of Locale Theory or in that of Topos Theory; we shall choose the latter since it is more general and represents the natural environment for investigating the relationships between the `generalized spaces' and their presentations, as well as for interpreting the topological and algebraic perspectives on the dualities in a unified way, that is in terms of topological and algebraic sites of definition for the same topos according to the philosophy `toposes as bridges' of \cite{OC10}.

\subsection{The topos-theoretic interpretation}\label{int}

Our topos-theoretic interpretation of Priestley duality stems from the observation that the essential part of the duality, namely the fact that the open sets of a coherent space can be recovered from the associated Priestley space as the open upper sets, can be naturally expressed in diagrammatic form, as follows.

The diagram
\[  
\xymatrix {
{\tau} \ar[d] \ar[rr] & & {\tau}_{pr}  \ar[d]   \\
Upp_{\leq_{X}}(X) \ar[rr] & & {\mathscr{P}}(X)}
\]
in the category of frames, where $X_{\tau}:=(X, {\tau})$ is a coherent space, $X_{{\tau}_{pr}}:=(X, {\tau}_{pr})$ is the associated Priestley space, $\leq_{X}$ is the specialization order on $X$ induced by the topology $\tau$, $Upp_{\leq_{X}}(X)$ is the frame of $\leq_{X}$-upper sets on $X$ and all the arrows are canonical inclusions, is a pullback. Indeed, by Priestley's duality (in the formulation given by Cornish in \cite{cornish}), the open sets in ${\tau}$ are precisely the sets lying in the intersection of ${\tau}_{pr}$ and $Upp_{\leq_{X}}(X)$ inside ${\mathscr{P}}(X)$. In fact, this diagram is the image, under the canonical functor $\textbf{Top}^{\textrm{op}}\to \textbf{Frm}$, of a pushout diagram in the category of topological spaces, obtained by regarding the frames appearing in the diagram as the frame of open sets of topological spaces whose underlying set is $X$.

The diagram above corresponds, via the equivalence of the category of locales with the category $\mathfrak{Loc}$ of localic toposes and geometric morphisms between them, to the following pushout in $\mathfrak{Loc}$:
\[  
\xymatrix {
[X, \Set] \ar[d]^{u} \ar[rr]^{\chi} & & \Sh(X_{{\tau}_{pr}}) \ar[d]^{f}   \\
[X_{\leq_{X}}, \Set] \ar[rr]^{\xi} & & \Sh(X_{\tau})}
\]
where 
\begin{enumerate}[(i)]
\item $X_{\tau}$ is the topological space whose underlying set is $X$ and whose frame of open sets is $\tau$,
\item $X_{{\tau}_{pr}}$ is the topological space whose underlying set is $X$ and whose frame of open sets is $\tau_{pr}$,
\item $X_{\leq_{X}}$ is the specialization order on $X$ induced by the topology $\tau$,
\item $f:\Sh(X_{{\tau}_{pr}}) \to \Sh(X_{{\tau}})$ is the geometric morphism induced by the continuous map of topological spaces $X_{{\tau}_{pr}}\to X_{{\tau}}$ whose underlying map is the identity on $X$, 
\item $\xi:[X_{\leq_{X}}, \Set] \to \Sh(X_{\tau})$ is the canonical geometric morphism induced by the specialization preorder on $X_{\tau}$,
\item $\chi:[X, \Set]\to \Sh(X_{{\tau}_{pr}})$ is the geometric morphism induced by the indexing of the set of points of the topos $\Sh(X_{{\tau}_{pr}})$ by the set $X$, and
\item $u:[X, \Set] \to [X_{\leq_{X}}, \Set]$ is the geometric morphism induced by the obvious functor $X\to X_{\leq_{X}}$ (where $X$ is considered a discrete category and $X_{\leq_{X}}$ is considered as a preorder category in the obvious way).
\end{enumerate} 
 
In the above diagram, the spaces $X_{{\tau}_{pr}}$ and $\Sh(X_{\tau})$ being sober, the set $X$ can be identified both with the set of points of the topos $\Sh(X_{{\tau}_{pr}})$ and with the set of points of the topos $\Sh(X_{\tau})$. In fact, this identification can be seen as being induced by composition of the points of the toposes $\Sh(X_{{\tau}_{pr}})$ and $\Sh(X_{\tau})$ with the geometric morphism $f$; specifically, the following diagram in $\Set$ commutes.
\[  
\xymatrix {
X \ar[d]^{\cong} \ar[rr]^{1_{X}} & & X \ar[d]^{\cong}   \\
Pts(\Sh(X_{\tau})) \ar[rr]^{h} & & Pts(\Sh(X_{{\tau}_{pr}})),}
\]
where $Pts(\Sh(X_{\tau}))$ (resp. $Pts(\Sh(X_{{\tau}_{pr}}))$) denotes the set of points of the topos $\Sh(X_{\tau})$ (resp. of the topos $\Sh(X_{{\tau}_{pr}})$), the map $h$ is the function between the points of the two toposes induced by composition with the morphism $f$, and the two isomorphisms are the canonical bijections identifying the points of a sober space with the points of the topos of sheaves on it.

\section{A general framework for Priestley-type dualities}\label{general}

The following definition represents a natural `invariant' generalization of the above setup. 

\begin{definition}
A \emph{Priestley context} consists of a geometric morphism $f:{\cal E}\to {\cal F}$ of localic toposes and two sets of points $X_{\cal E}$ and $X_{\cal F}$ respectively of the topos $\cal E$ and of the topos $\cal F$ such that composition with $f$ induces a bijection $X_{f}:X_{{\cal E}}\to X_{{\cal F}}$ with the property that the diagram    
\[  
\xymatrix {
[X_{\cal E}, \Set] \ar[d]^{u} \ar[rr]^{\chi} & & {\cal E}  \ar[d]^{f}   \\
[{X_{\cal E}}_{\leq_{\cal F}}, \Set] \ar[rr]^{\xi} & & {\cal F},}
\]
where
\begin{enumerate}[(i)]
\item ${X_{\cal E}}_{\leq_{\cal F}}$ is the preorder on ${X_{\cal E}}$ corresponding to the specialization preorder ${X_{\cal F}}_{\leq_{\cal F}}$ on ${X_{\cal E}}$ under the bijection $X_{f}$,
\item the geometric morphism $\xi$ is given by the composite of the equivalence $[{X_{\cal E}}_{\leq_{\cal F}}, \Set]\to [{X_{\cal F}}_{\leq_{\cal F}}, \Set]$ induced by the isomorphism $X_{f}:{X_{\cal E}}_{\leq_{\cal F}} \simeq {X_{\cal F}}_{\leq_{\cal F}}$ with the canonical morphism $[{X_{\cal F}}_{\leq_{\cal F}}, \Set] \to {\cal F}$,
\item $\chi:[X_{\cal E}, \Set] \to {\cal E}$ is the geometric morphism corresponding to the set of points $X_{\cal E}$ of the topos $\cal E$ and
\item $u:[X_{\cal E}, \Set] \to [{X_{\cal E}}_{\leq_{\cal F}}, \Set]$ is the geometric morphism induced by the canonical functor $X_{\cal E} \to {X_{\cal E}}_{\leq_{\cal F}}$, 
\end{enumerate}
is a pushout in the category $\mathfrak{Loc}$.
\end{definition}

We remarked in section \ref{int} that the classical Priestley duality gives rise to a Priestley context for each distributive lattice. In this section, we show that one can easily generate other instances of Priestley contexts involving structures other than distributive lattices, and that in fact there is a uniform way for building them.

\subsection{Generalized patch topologies}

Our general method for building Priestley contexts is based on a generalization of the concept of patch topology.

\begin{definition}
Let $A$ be a subset of a powerset ${\mathscr{P}}(X)$. The $A$-subbasic topology on $X$ is the topology on $X$ having $A$ as a subbasis. The $A$-patch topology on $X$ is the topology on $X$ having as a sub-basis the set consisting of the elements in $A$ and their complements in ${\mathscr{P}}(X)$.
\end{definition}

We shall denote the topological space obtained by endowing $X$ with the $A$-subbasic topology (resp. with the $A$-patch topology) by $X_{A}$ (resp. by $X_{A}^{p}$). Notice that 

\begin{remarks}\label{funpr}
\begin{enumerate}[(a)]

\item The subsets of $X$ which can be expressed as finite unions of intersections of elements of ${\mathscr{P}}(X)$ which either belong to $A$ or are a complement of a subset belonging to $A$ form a basis of clopen subsets for the topology $X_{A}^{p}$. We shall refer to this collection of subsets, endowed with the natural subset-inclusion ordering, as the \emph{Priestley Boolean algebra} generated by $(A, X)$. 

\item The construction of the patch topology can be made functorial. Suppose that $A$ is a subset of ${\mathscr{P}}(X)$, $B$ is a subset of ${\mathscr{P}}(Y)$, and $f:X\to Y$ is a function such that $f^{-1}:{\mathscr{P}}(Y)\to {\mathscr{P}}(X)$ restricts to a map ${\cal O}(Y_{B})\to {\cal O}(X_{A})$. Then $f^{-1}$ restricts to a map ${\cal O}(Y_{B}^{p})\to {\cal O}(X_{A}^{p})$. In particular, every continuous map of subbasic topologies extends to a continuous map of the associated patch topologies.
 
\item Let ${\cal C}\subseteq {\mathscr{P}}(X)$ and $\cal L$ be the category $\textbf{Bool}$ whose objects are the Boolean algebras and whose arrows are the Boolean algebra homomorphisms between them. Then the $\cal L$-substructure $G^{X}_{({\cal C}, {\cal L})}$ of ${\mathscr{P}}(X)$ generated by $\cal C$ coincides with the Priestley Boolean algebra generated by $({\cal C}, X)$.  
\end{enumerate}
\end{remarks}

\begin{lemma}\label{lemmapr}
Let $A$ be a subset of a powerset ${\mathscr{P}}(X)$. If the space $X_{A}^{p}$ is compact then the open sets of the space $X_{A}$ are exactly the open sets of the space $X_{A}^{p}$ which are upper-sets with respect to the specialization preorder on $X$ corresponding to the topological space $X_{A}$.
\end{lemma}

\begin{proofs}
We can show our thesis by generalizing the well-known argument for the classical Priestley duality (cf. for example Proposition II4.6 \cite{stone}). Clearly, every open set of $X_{A}$ is both an open set of $X_{A}^{p}$ and an upper-set with respect to the specialization preorder $\leq_{A}$ on $X$ induced by $X_{A}$, so it remains to prove the converse direction. Let $U$ be an open set of $X_{A}^{p}$ which is a $\leq_{A}$-upper set; we will show that $U$ is an open set of $X_{A}$ by providing, for every fixed point $x$ of $U$, an open set of $X_{A}$ containing $x$ and contained in $U$. Let $x$ be an element of $U$; then for each $y\in X\setminus U$, we have that $x\nleq y$, from which it follows that there exists an open set $V_{y}$ of $X_{A}$ belonging to $A$ (notice that, by our hypotheses, $A$ is a subbasis for $X_{A}$) such that $x\in V_{y}$ and $y\notin V_{y}$ (equivalently, $y\in X\setminus V_{y}$). The sets $\{X\setminus V_{y} \textrm{ | } y\in X\setminus U \}$ are thus open sets of $X_{A}^{p}$ which jointly cover $X\setminus U$; from the compactness of $X_{A}^{p}$ it thus follows that this covering family admits a finite subcover  $\{X\setminus V_{y_{i}} \textrm{ | } i=1, \ldots, n \}$, equivalently the intersection of the $V_{y_{i}}$ is an open neighborhood of $x$ contained in $U$. 
\end{proofs}

Under natural hypotheses, we can characterize the topological spaces arising by putting the $A$-patch topology on a set $X$ in a natural way, using the specialization order $\leq_{A}$ induced by the topology $X_{A}$. Specifically, we have the following result.

\begin{theorem}\label{prchar}
Let $A$ be a subset a powerset ${\mathscr{P}}(X)$, and $(X, \tau, \leq)$ be a preordered compact topological space. Then $\tau$ is the $A$-patch topology on $X$ and $\leq$ coincides with the specialization preorder $\leq_{A}$ induced by the topology $X_{A}$ if and only if every set in $A$ is a $\leq$-upper set and $\tau$ satisfies the following separation axiom: for any $x, y\in X$ such that $x\nleq_{A} y$, there is a clopen $\leq$-upper set $U$ of $\tau$ belonging to $A$ such that $x\in U$ and $y\notin U$. 
\end{theorem} 

\begin{proofs}
The `only if' direction is obvious from the definition of $A$-patch topology on $X$, so it remains to prove the converse one. Suppose that $\tau$ satisfies the separation axiom in the statement of the theorem. Given an open set $U$ of $(X, \tau)$ and a point $x$ of $U$, for each $y\in X\setminus U$, either $y\nleq x$ or $x\nleq y$. Then there exists a clopen set of $\tau$ in $A$ which contains $x$ and misses $y$ or a clopen set of $\tau$ given by the complement of a set in $A$ which contains $x$ and misses $y$. Obviously the intersection of these clopen neighborhoods of $x$ (for $y$ varying in $X\setminus U$) does not meet $X\setminus U$. Therefore, as $(X, \tau)$ is compact, there exists a finite intersection of these clopen neighborhoods of $x$ missing $X\setminus U$. This finite intersection is a clopen neighborhood $C$ of $x$ contained in $U$ which is a finite intersection of sets in $A$ and their complements. This proves that $\tau$ is the $A$-patch topology on $X$. 

To show that $\leq$ coincides with the specialization preorder $\leq_{A}$ we observe that, by our hypotheses, $(X, \tau, \leq)$ is compact and satisfies the Priestley separation axiom, and hence it is a Priestley space. Therefore, in view of the classical Priestley duality, it suffices to prove that the $\leq$-upper open sets in $\tau$ coincide precisely with the open sets in $X_{A}$. Given a $\leq$-upper open set $U$ of $(X, \tau)$ and a point $x$ of $U$, for each $y\in X\setminus U$, we have $y\nleq x$ (since $U$ is a $\leq$-upper set); then, by the separation axiom in the statement of the theorem, there exists a clopen $\leq$-upper set of $\tau$ in $A$ which contains $x$ and misses $y$ and hence, as $(X, \tau)$ is compact, there exists a finite intersection of these clopen neighborhoods of $x$ missing $X\setminus U$. This finite intersection is a clopen neighborhood $C$ of $x$ contained in $U$ which is a finite intersection of sets in $A$ and hence an open set of $X_{A}$; thus $U$ is an open set of $X_{A}$, since for any point $x$ of $U$ there exists a neighborhood of $x$ contained in $U$ and belonging to $X_{A}$. This argument shows that any $\leq$-upper open set in $\tau$ belongs to $X_{A}$; notice in passing that if $U$ is a ($\leq$-upper) clopen then $U$ is a \emph{finite} union of finite intersections of sets in $A$. The fact that every open set in $X_{A}$ is a $\leq$-upper set in $\tau$ follows from the fact that, by our hypothesis, every set in $A$ is a $\leq$-upper set. 
\end{proofs}

\begin{corollary}\label{generationPriestley}
Let $A$ be a subset a powerset ${\mathscr{P}}(X)$. Then, provided that it is compact, the preordered topological space $(X_{A}^{p}, \leq_{A})$ is a Priestley space, whose corresponding distributive lattice is the sublattice of $\mathscr{P}(X)$ consisting of the subsets which are finite unions of finite intersections of subsets in $A$ and whose corresponding coherent space is $X_{A}$.
\end{corollary}

\begin{proofs}
The corollary follows immediately from the arguments in the proof of Theorem \ref{prchar}.
\end{proofs}

\subsection{Free structures and Priestley contexts}

In this section we shall see that free structures provide a very natural way for building Priestley contexts.

\begin{theorem}\label{patch}
Let $\cal C$, $\cal M$, ${\cal L}=\textbf{Bool}$, $m:{\cal C}\to L$, $X_{\cal C}$ $J_{\cal C}$, $Y_{L}$ be a set of data satisfying the hypotheses of Theorem \ref{prgen}. Then, if we regard $\cal C$ as a subset of $\mathscr{P}(X_{\cal C})$ via the composite map ${\cal C}\mono Id_{J_{\cal C}}({\cal C})\cong {\cal O}(X_{\cal C})\subseteq \mathscr{P}(X_{\cal C})$, there is a set of points $Z_{L}$ of the topos $\Sh(X_{{\cal C}}^{p})$ such that the geometric morphism $\Sh(X_{{\cal C}}^{p})\to \Sh(X_{{\cal C}})$ induced by the continuous function $X_{{\cal C}}^{p} \rightarrow X_{{\cal C}}$ whose underlying map is the identity on $X$, together with the set of points $X_{\cal C}$ of the topos $\Sh(X_{{\cal C}})$ and the set of points $Z_{L}$ of the topos $\Sh(X_{{\cal C}}^{p})$ is a Priestley context. In particular, if the class $\cal M$ is axiomatized by the theory of $J_{\cal C}$-prime filters on $\cal C$ then $Z_{L}$ can be taken to be the set of all the points of the topos $\Sh(X_{{\cal C}}^{p})$.
\end{theorem}

\begin{proofs}
Below we shall employ the notions and definitions used in Theorem \ref{prgen} and its proof.
 
Clearly, the topological space obtained by equipping the given set of points of the topos $\Sh({\cal C}, J_{\cal C})$ with the subterminal topology coincides with the space obtained by equipping $X_{\cal C}$ with the $\cal C$-subbasic topology, where $\cal C$ is regarded as a subset of $\mathscr{P}(X_{\cal C})$ via the composite map ${\cal C}\mono Id_{J_{\cal C}}({\cal C})\cong {\cal O}(X_{\cal C})\subseteq \mathscr{P}(X_{\cal C})$; in particular, since the topos $\Sh({\cal C}, J_{\cal C})$ has enough points, we have an equivalence $\Sh({\cal C}, J_{\cal C})\simeq \Sh(X_{\cal C})$. 

On the other hand, from the proof of Theorem \ref{prgen} we know that there is a commutative diagram
\[  
\xymatrix {
G^{X_{\cal C}}_{({\cal C}, {\cal L})}  \ar[d]^{\cong} \ar[rr] & & \mathscr{P}(X_{\cal C}) \ar[d]^{\mathscr{P}(u)}   \\
L \ar[rr]^{k} & & {\mathscr{P}}(Y_{L}),}
\]
where the arrows $u$ and $k$ are those defined in the proof of the theorem. From this it follows that the space $Z$ obtained by endowing the set of points $Y_{L}$ of the topos $\Sh(L, K_{L})$ with the subterminal topology is homeomorphic to the space $X_{{\cal C}}^{p}$, since the bijection $u:X_{\cal C}\cong Y_{L}$ between the points of $X_{{\cal C}}^{p}$ and the points of the space $Z$ maps basic open sets of $X_{{\cal C}}^{p}$ (that is, subsets in $G^{X_{\cal C}}_{({\cal C}, {\cal L})}$, cf. Remark \ref{funpr}(c)) to basic open sets of $Z$ (that is, subsets of the form $k(l)$ for an element $l\in L$); in particular, since the topos $\Sh(L, K_{L})$ has enough points, we have an equivalence $\Sh(L, K_{L})\simeq \Sh(X_{{\cal C}}^{p})$. Let us set $Z_{L}$ equal to the set of points of the topos $\Sh(X_{{\cal C}}^{p})$ corresponding to the set of points $Y_{L}$ of the topos $\Sh(L, K_{L})$ under this equivalence.  
The fact that the geometric morphism $\Sh(X_{{\cal C}}^{p})\to \Sh(X_{{\cal C}})$ induced by the continuous function $X_{{\cal C}}^{p} \rightarrow X_{{\cal C}}$ whose underlying map is the identity on $X$ is, together with the set of points $X_{{\cal C}}$ and $Z_{L}$ respectively of the topos $\Sh(X_{{\cal C}})$ and of the topos $\Sh(X_{{\cal C}}^{p})$ is a Priestley context thus follows from Lemma \ref{lemmapr} by arguing as in section \ref{int}.     
The last part of the theorem follows from the obvious remark that if $X_{\cal C}$ is the set of \emph{all} the points of the topos $\Sh({\cal C}, J_{\cal C})$ then, $L$ being the free Boolean algebra on $\cal C$, $Y_{L}$ is the set of all the points of the topos $\Sh(L, K_{L})$.
\end{proofs}

We notice that if $K_{L}$ is equal to the coherent topology on $L$ and $Y_{L}$ is equal to the set of all the points of the topos $\Sh(L, K_{L})$ then the space $X_{{\cal C}}^{p}$ is a Stone space (since it is homeomorphic to the space of points of the topos $\Sh(L, K_{L})$). In fact, in order to build real analogues of Priestley-type dualities involving \emph{Stone} spaces, we shall take $K_{L}$ equal to the coherent topology on $L$, and $X_{\cal C}$ to consist of \emph{all} the points of the topos $\Sh({\cal C}, J_{\cal C})$.

\subsection{The general method}\label{genmeth}

Let $\cal K$ be a category of posets, each of which equipped with a subcanonical Grothendieck topology such that the morphisms of the category $\cal K$ coincide with the morphisms of the associated sites. We denote by $J_{\cal C}$ the Grothendieck topology associated to a poset $\cal C$ in $\cal K$. 

Let us moreover suppose that all the theories of $J_{\cal C}$-prime filters on $\cal C$ introduced in \cite{OC11} are coherent (notice that this is always the case if $\cal C$ is a meet-semilattice), so that all the toposes $\Sh({\cal C}, J_{\cal C})$ have enough points; then we have, by the results of \cite{OC11}, a functor $A:{\cal K}^{\textrm{op}}\to \textbf{Top}$, assigning to any structure $\cal C$ in $\cal K$ the space of points of the topos $\Sh({\cal C}, J_{\cal C})$. By Remark \ref{funpr}(b), the generalized patch topology construction allows us to lift this functor to one with values in the category $\textbf{PTop}$ of preordered topological spaces. Specifically, for any $\cal C$ in $\cal K$, denoted by $X_{\cal C}$ the space of points of the topos $\Sh({\cal C}, J_{\cal C})$, and regarded $\cal C$ as a subset of $\mathscr{P}(X_{\cal C})$ via the embedding given by the composite map ${\cal C}\mono Id_{J_{\cal C}}({\cal C})\cong {\cal O}(X_{\cal C})\subseteq \mathscr{P}(X_{\cal C})$, we set $A_{p}({\cal C})$ equal to the set $X_{\cal C}$, equipped with $\cal C$-patch topology and the specialization preorder $\leq_{\cal C}$ corresponding to the space $X_{\cal C}$, and for any arrow $f$ in $\cal K$ we set $A_{p}(f)$ equal to the map $A(f)$, regarded as a morphism of preordered topological spaces.

Recall from section \ref{pts} that we have a functor $r_{P}:\textbf{PTop}\to\textbf{Top}$ which is left adjoint to the canonical inclusion functor $i_{P}:\textbf{Top}\to\textbf{PTop}$ and is such that the composite $r_{P}\circ i_{P}$ is isomorphic to the identity functor on $\textbf{Top}$. In fact, the assignments above give rise to a functor $A_{p}:{\cal K}^{\textrm{op}}\to \textbf{PTop}$ such that $r_{P}\circ A_{p}\cong A$:

\[  
\xymatrix {
 & & \textbf{PTop} \ar[d]^{r_{P}}   \\
{\cal K}^{\textrm{op}} \ar[urr]^{A_{p}} \ar[rr]^{A} & & \textbf{Top}}
\]

Recall from \cite{OC11} that if $F:{\cal U}\to {\cal V}$ is a functor which creates isomorphisms then there exists a smallest subcategory $\cal W$ of $\cal V$ closed under isomorphisms in $\cal V$ such that $F$ factors through the canonical inclusion functor ${\cal W} \hookrightarrow {\cal V}$. Such category is called the extended image of the functor $F$ and denoted by $ExtIm(F)$.

Notice that if $A$ is faithful then $A_{p}$ is faithful as well. Also, since the functor $r_{P}$ is injective on arrows, if $A$ creates isomorphisms then $A_{p}$ creates isomorphisms as well; in particular, if $A$ is part of an equivalence with its extended image then $A_{p}$ is also part of an equivalence with its extended image; indeed, the commutativity of the above triangle forces $r_{P}$ to send the extended image of $A_{p}$ to the extended image of $A$ and hence a categorical left (and hence right as well, by Proposition 3.4 \cite{OC11}) inverse for $A_{p}$, defined on its extended image is given by the composite of the categorical inverse of $A$ with $r_{P}$.    

Suppose that for any $\cal C$ in $\cal K$ we have a set of maps ${\cal M}_{\cal C}$ from $\cal C$ to Boolean algebras satisfying the property that for any function from $\cal C$ to a Boolean algebra if the composition of it with each of the maps belonging to a jointly injective family of Boolean algebra homomorphisms belongs to ${\cal M}_{\cal C}$ then the function itself belongs to ${\cal M}_{\cal C}$, and that for any $\cal C$ in $\cal K$ the arrows ${\cal C}\to \{0,1\}$ in ${\cal M}_{\cal C}$ coincide precisely with the $J_{\cal C}$-continuous flat functors ${\cal C}\to \{0,1\}$. Then, if for each $\cal C$ in $\cal K$ the free $(\textbf{Bool}, {\cal M}_{\cal C})$-structure $L_{\cal C}$ on $\cal C$ exists, the hypotheses of Theorem \ref{patch} are satisfied, with $\cal L$ being the category $\textbf{Bool}$ of Boolean algebras and $K_{L}$ being the coherent topology on $L$ for each Boolean algebra $L$. In particular, under these assumptions it follows from the proof of theorem \ref{patch} that there is an equivalence of toposes $\Sh(X^{p}_{\cal C})\simeq \Sh(L_{\cal C}, K_{L_{\cal C}})$, from which it follows that the space $X^{p}_{\cal C}$ is compact. Since from the proof of Theorem \ref{prchar} we know that the spaces $(X^{p}_{\cal C}, \leq_{\cal C})$ satisfy the Priestley separation axiom, we can conclude that under the assumptions specified above all the preordered topological spaces in the image of the functor $A_{p}$ are Priestley spaces. 
   
If the Grothendieck topologies are $C$-induced for a topos-theoretic invariant $C$ satisfying the hypotheses of Theorem 3.25 \cite{OC11} then the functor $A:{\cal K}^{\textrm{op}} \to \textbf{Top}$ yields a `Stone-type' duality between $\cal K$ and a subcategory $ExtIm(A)$ of $\textbf{Top}$, and hence the functor $A_{p}:{\cal K}^{\textrm{op}} \to \textbf{PTop}$ yields a `Priestley-type' duality between $\cal K$ and a category of Priestley spaces, given by $ExtIm(A_{p})$. 

Let us now consider the problem of characterizing the subcategories of $\textbf{PTop}$ arising as the extended images of the functors $A_{p}$. Notice that, if all the spaces in $ExtIm(A)$ are coherent then, by the classical Priestley duality, the Priestley spaces in $ExtIm(A_{p})$ can be characterized as the Priestley spaces $(X, \tau, \leq)$ such that $r_{P}(X, \tau, \leq)$ belongs to $ExtIm(A)$, while the arrows in $ExtIm(A_{p})$ can be characterized as the arrows whose image under $r_{P}$ lies in the extended image $ExtIm(A)$. More explicitly, by Theorem \ref{prchar}, if the spaces in $ExtIm(A)$ can be characterized as the sober topological spaces with a basis of $C$-compact open sets satisfying some property $P$, the spaces in $ExtIm(A_{p})$ can be characterized as the compact ordered topological spaces which satisfy the following separation axiom: for any $x, y\in X$ such that $x\nleq y$, there is a clopen $\leq$-upper set $U$ of $\tau$ which is $C$-compact among the $\leq$-upper open sets of $\tau$ such that $x\in U$ and $y\notin U$ and which moreover satisfy the property that the collection of their $\leq$-upper clopen sets which are $C$-compact among the $\leq$-upper clopen sets satisfy property $P$. Alternatively, they can be characterized as the sober ordered topological spaces such the subsets of their underlying sets which are $\leq$-upper clopen and $C$-compact among the $\leq$-upper open sets form, together with their complements, a subbasis of the space satisfying property $P$. The arrows in $ExtIm(A_{p})$ can be characterized as the arrows in $\textbf{PTop}$ such that the inverse image of any set which is $\leq$-upper clopen and $C$-compact among the $\leq$-upper open sets is $C$-compact in such a way that the restriction of the inverse image function to such subsets can be identified with an arrow in $\cal K$ (notice that this latter condition can be dropped in the case the morphisms ${\cal C}\to {\cal D}$ in $\cal K$ coincide with the morphisms $({\cal C}, J_{\cal C})\to ({\cal D}, J_{\cal D})$ of the associated sites).      

We shall see concrete examples of `Priestley-type' dualities generated through the method described above in section \ref{ex}. 
 
\subsubsection{The algebraic interpretation}\label{algint}

If we identify ordered Stone spaces with Boolean algebras equipped with an order on the points of their spectra according to Stone duality for Boolean algebras, we can define the functor $A_{p}$, and characterize its extended image, in algebraic terms. 

Specifically, we can think of a Priestley-type space $(X, \tau, \leq)$ as a pair $(B, \leq)$, where $\leq$ is an ordering on the Stone spectrum $Spec(B)$ of $B$ satisfying the Priestley separation axiom. Notice that if $X=Spec(B)$ then the points $x$ of $X$ can be identified with the prime filters $F$ on $P$ and hence the condition that for any $x,x'\in X$ such that $x\nleq x'$ there should exist a clopen $\leq$-upper set $U$ such that $x\in U$ and $x'\notin U$ rewrites as follows: for any $F, F'\in Spec(B)$ such that $F\nleq F'$ there exists an element $b\in B$ such that $b\in F$, $b\notin F'$ and for any prime filters $G,G'\in Spec(B)$ with the property that $G\leq G'$, $b\in G$ implies $b\in G'$. We shall say that an element $b\in B$ is \emph{$\leq$-upper} if it satisfies the condition that for any prime filters $G,G'\in Spec(B)$ with the property that $G\leq G'$, $b\in G$ implies $b\in G'$. 

This remark, combined with the algebraic construction of Priestley spaces through free structures, paves the way for an entirely algebraic reformulation of the Priestley-type dualities obtained through the method of section \ref{genmeth} above.

Let $\cal K$ be a category of posets $\cal C$ satisfying the hypotheses of the method of section \ref{genmeth}. For $\cal C$ in $\cal K$, let us denote by $B_{\cal C}$ the free $({\cal M}_{\cal C}, \textbf{Bool})$-structure on $\cal C$. Let us denote by $i_{{\cal C}}:{\cal C}\to B_{\cal C}$ the universal map from $\cal C$ to the free Boolean algebra $B_{\cal C}$. Clearly, by the universal property of the algebras $B_{\cal C}$, any morphism of sites $({\cal C}, J_{\cal C})\to ({\cal C}', J_{{\cal C}'})$ (that is, any morphism $f:{\cal C}\to {\cal C}'$ in $\cal K$) induces a Boolean algebra homomorphism $B_{f}:B_{\cal C}\to B_{{\cal C}'}$, which in fact is the unique Boolean algebra homomorphism $s:B_{\cal C}\to B_{{\cal C}'}$ such that $s\circ i_{\cal C}=i_{{\cal C}'}\circ f$. 

We thus have a functor
\[
B:{\cal K}\to \textbf{Bool}_{\leq},
\]  
where $\textbf{Bool}_{\leq}$ is the category whose objects are the pairs $(B, \leq)$, where $B$ is a Boolean algebra and $\leq$ is a order on the set $Spec(B)$ of prime filters of $B$ with the property that for any $F, F'\in Spec(B)$ such that $F\nleq F'$ there exists a $\leq$-upper element $b\in B$ such that $b\in F$ and $b\notin F'$, and whose arrows $(B, \leq)\to (B', \leq')$ are the Boolean algebra homomorphisms $f:B\to B'$ such that $f^{-1}:Spec(B')\to Spec(B)$ is order-preserving (i.e., for any $F,F'\in Spec(B)$, $F\leq F'$ in $Spec(B)$ implies $f^{-1}(F)\leq' f^{-1}(F')$ in $Spec(B')$). 

The functor $B$ is defined as follows: for any $\cal C$ in $\cal K$, $B({\cal C})=(B_{\cal C}, \leq_{\cal C})$, where $\leq_{\cal C}$ is the order on $Spec(B)$ given by: for any $F, F'\in Spec(B)$, $F\leq_{\cal C} F'$ if and only if $F\cap {\cal C}\subseteq F'\cap {\cal C}$, while for any arrow $f:{\cal C}\to {\cal C}'$ in $\cal K$, $B(f)=B_{f}$.

Recall from \cite{OC11} (specifically, Definition 3.15) that, given a frame-theoretic invariant property $C$ of families of elements of a frame, an element of a frame $L$ is said to be \emph{$C$-compact} if every covering family of $l$ in $L$ has a refinement satisfying $C$; that is, whenever $a=\mathbin{\mathop{\textrm{\huge $\vee$}}\limits_{i\in I}}a_{i}$ in $L$ there exists a family $\{b_{j} \leq a \textrm{ | } j\in J\}$ of elements of $L$ satisfying $C$ such that for every $j\in J$ there exists $i\in I$ such that $b_{j}\leq a_{i}$, and the join $\mathbin{\mathop{\textrm{\huge $\vee$}}\limits_{i\in I}}a_{i}$ in $L$ is equal to the join $\mathbin{\mathop{\textrm{\huge $\vee$}}\limits_{j\in J}}b_{j}$ in $L$. We shall apply this notion to the frame of open sets of the coherent spaces associated to the Priestley spaces under consideration.
 
Notice that, given a Priestley space $(X, \tau, \leq)$, an $\leq$-upper open set $U$ in $\tau$ is $C$-compact among the $\leq$-upper open sets (that is, it is $C$-compact in the coherent space associated to the Priestley space) if and only if every covering in $U$ by clopen $\leq$-upper sets has refinement satisfying $C$ (since the clopen $\leq$-upper sets form a basis for the coherent space associated to the Priestley space); in fact, we can equivalently require the refinement to consist of clopen $\leq$-upper sets, since the topologies $J_{\cal C}$ satisfy by our the hypotheses the condition that if a covering family consisting of principal $J_{\cal C}$-ideals on $\cal C$ admits a refinement satisfying $C$ then it also admits a refinement satisfying $C$ and consisting of principal $J_{\cal C}$-ideals on $\cal C$, and these ideals are sent by the embedding $Id_{J_{\cal C}}({\cal C})\mono {\cal O}(X_{\cal C})\hookrightarrow Spec(B_{\cal C})$ to clopen $\leq$-upper sets of $Spec(B_{\cal C})$, where $B({\cal C})=(B_{\cal C}, \leq)$. Therefore, an element $b\in B$ is $\leq$-upper in $Spec(B)$ and $C$-compact among the $\leq$-upper sets of $Spec(B)$ if and only if for any prime filters $G,G'\in Spec(B)$ with the property that $G\leq G'$, $b\in G$ implies $b\in G'$ and for any finite set of $\leq$-upper elements $\{b_{1}, \ldots, b_{n}\}$ such that $b=b_{1}\vee \cdots \vee b_{n}$ in $B$ there exists a refinement of the family $\{b_{1}, \ldots, b_{n}\}$ satisfying the invariant $C$. Indeed, by the compactness of $b$ in $Spec(B)$, any covering of $b$ in $Spec(B)$ admits a finite subcovering by elements of $B$, and a family $\{b_{1}\leq b, \ldots, b_{n}\leq b\}$ covers $b$ in $Spec(B)$ if and only if $b=b_{1}\vee \cdots \vee b_{n}$ in $B$. For this reason, we shall call the elements $b\in B$ which are $\leq$-upper in $Spec(B)$ and $C$-compact the finitely $C$-compact elements among the $\leq$-upper elements of $B$. 

As we observed in section \ref{genmeth} above, if the Grothendieck topologies are $C$-induced for a topos-theoretic invariant $C$ satisfying the hypotheses of Theorem 3.25 \cite{OC11} then the functor $B$ defined above yields an equivalence of categories onto its extended image. In fact, any $\cal C$ in $\cal K$ can be recovered from the Boolean algebra $B({\cal C})$ as the set of its elements which are $\leq$-upper and finitely $C$-compact among the $\leq$-upper elements of $B$, while any arrow $f:{\cal C}\to {\cal C}'$ in $\cal K$ can be recovered as the restriction of the corresponding arrow $B(f)$ along the canonical inclusions $i_{\cal C}$ and $i_{{\cal C}'}$.  
    
Let us now turn to the problem of characterizing the extended image of the functor $B$. The arrows in $ExtIm(B)$ can be characterized as the arrows in $\textbf{Bool}_{\leq}$ between objects in $ExtIm(B)$ which send any $\leq$-upper element $b\in B$ with the property that every finite covering of it by $\leq$-upper elements admits a refinement by a family which satisfies $C$ to a $\leq$-upper element of $B'$ with the property that every finite covering of it by $\leq$-upper elements admits a refinement by a family which satisfies $C$. The following theorem takes care of characterizing the objects in $ExtIm(B)$.

\begin{theorem}\label{tfae}
Under the hypotheses above, the following conditions are equivalent.
\begin{enumerate}[(i)]
\item $(B, \leq)$ is an object of $ExtIm(B)$;

\item For any $F, F'\in Spec(B)$ such that $F\nleq F'$ there exists a $\leq$-upper element $b\in B$ with the property that $b\in F$ and $b\notin F'$ and every finite covering of it in $B$ by $\leq$-upper elements admits a refinement by a family which satisfies $C$, and the set of $\leq$-upper elements of $B$, with the induced order, can be identified with a structure in $\cal K$;

\item The inclusion $B^{\ast}\subseteq B$ where $B^{\ast}$ is the subset of $B$ consisting of the $\leq$-upper elements $b$ of $B$ which satisfy the property that every finite covering of it by $\leq$-upper elements admits a refinement by a family satisfying $C$ realizes $B$ as the free $({\cal M}, \textbf{Bool})$-algebra on $B^{\ast}$ and satisfies the property that for any $F, F'\in Spec(B)$ such that $F\nleq F'$ there exists an element $b\in B^{\ast}$ with the property that $b\in F$ and $b\notin F'$. 
\end{enumerate}
\end{theorem} 

\begin{proofs}
The equivalence $(i)\biimp (ii)$ represents the algebraic formulation of the topological characterization obtained in section \ref{genmeth} above.

The implication $(i)\imp (iii)$ follows from the characterization of the algebra $B_{\cal C}$ as the free $({\cal M}, \textbf{Bool})$-algebra on $\cal C$ and from the isomorphism ${B_{\cal C}}^{\ast}\cong {\cal C}$. To show that $(iii)$ implies $(i)$, we observe that if the inclusion $B^{\ast}\subseteq B$ realizes $B$ as the free $({\cal M}, \textbf{Bool})$-algebra on $B^{\ast}$ then $(B, \leq)=B({B_{\cal C}}^{\ast})$. Indeed, since $B$ is the free $({\cal M}, \textbf{Bool})$-algebra on $B^{\ast}$ then $B_{{B_{\cal C}}^{\ast}}=B$, while the order $\leq$ coincides with the order $\leq_{{B_{\cal C}}^{\ast}}$ since for any prime filters $F, F'\in Spec(B)$, $F\leq F'$ if and only if there is not an element $b\in B^{\ast}$ such that $b\in F$ and $b\notin F'$, equivalently if $F\cap B^{\ast} \subseteq F'\cap B^{\ast}$, i.e. if and only if $F\leq_{B^{\ast}} F'$.     
\end{proofs}

The following corollary provides, given a structure $\cal C$ and a class of maps $\cal M$ from $\cal C$ to Boolean algebras as above, a criterion for an inclusion $i:{\cal C}\hookrightarrow B$ of $\cal C$ into a Boolean algebra $B$ to realize $B$ as the free $({\cal M}, \textbf{Bool})$-algebra on $\cal C$.   
\begin{corollary}\label{corollarycar}
Under the hypotheses above, let $i:{\cal C}\hookrightarrow B$ be an inclusion, and $\leq_{\cal C}$ be the order on $Spec(B)$ defined above. Then $i$ realizes $B$ as the free $({\cal M}, \textbf{Bool})$-algebra on $\cal C$ if and only if $\cal C$ can be identified, via $i$, with the subset of $B$ consisting of the $\leq_{\cal C}$-upper elements $b$ of $B$ which satisfy the property that every finite covering of it by $\leq_{\cal C}$-upper elements admits a refinement by a family which satisfies $C$.
\end{corollary}
     
\begin{proofs}
If $i$ realizes $B$ as the free $({\cal M}, \textbf{Bool})$-algebra on $\cal C$ then, by the algebraic version of our method for building Priestley-type dualities established above, $\cal C$ can be identified with $B^{\ast}$. Conversely, if the condition in the corollary holds then condition $(ii)$ of Theorem \ref{tfae} holds for the space $(B, \leq_{\cal C})$, whence we conclude, by condition $(iii)$ in the theorem, that $i$ realizes $B$ as the free $({\cal M}, \textbf{Bool})$-algebra on $\cal C$.     
\end{proofs}

\subsection{Free structures and dualities}

In view of Theorem \ref{prgen}, free structures provide a natural way for building Priestley-type dualities. In fact, the idea of (functorially) recovering a given poset structure in a category $\cal K$ from a free structure on it (possibly equipped with some additional information such as for example an ordering on the set of points of its spectrum), and to characterize in intrinsic terms the structures which are free on the structures in $\cal K$ is very general and has applications well beyond the context specifically addressed in this paper. Also, it is natural to expect such dualities to admit extensions to reflections or coreflections to larger categories on which it is possible to define an assignment with values in the source category which naturally extends the operation of recovering a structure from the corresponding free structure on it. For example, in the classical Priestley duality one recovers a distributive lattice from the free Boolean algebra on it equipped with a partial order on its spectrum (as the set of its upper elements), while in the classical Stone duality one recovers a distributive lattice from the free frame on it (as the set of its compact elements). In fact, all the Priestley-type dualities established above in this paper arise, if viewed algebraically, from the construction of free Boolean algebras on particular kinds of posets, while the Stone-type dualities established in \cite{OC11} arise from the construction of free frames or posets on given preordered structures.

In this section we give further illustrations of this general phenomenon, by establishing new categorical equivalences defined by free functors. With time, new examples of dualities arising from free structures will likely arise and further enrich the picture.  

Let us work under the hypotheses of the last part of Theorem \ref{prgen}. In such a context, we have an equivalence of toposes
\[
\Sh(L, K_{L}) \simeq \Sh({\cal C}, J_{\cal C}),
\]  
induced by the canonical map $m:{\cal C}\to L$, regarded as a morphism of sites $({\cal C}, J_{\cal C})\to (L, K_{L})$. Notice that, both $J_{\cal C}$ and $K_{L}$ being subcanonical, $J_{\cal C}$ can be identified with the Grothendieck topology induced by $K_{L}$ on $\cal C$ under the canonical embedding ${\cal C}\hookrightarrow L$. So, by the results of \cite{OC11}, if the topology $J_{\cal C}$ is $C$-induced (in the sense of \cite{OC11}), $\cal C$ can be recovered from $L$ as the set of $C$-compact subterminals of the topos $\Sh(L, K_{L})$, and for any structure $L$ in $\cal L$ such that the set of its $C$-compact elements (i.e. the set of the $C$-compact subterminals on the topos $\Sh(L, K_{L})$ which are of the form $y(l)$ for an element $l\in L$, where $y$ is the Yoneda embedding $L\to \Sh(L, K_{L})$) can be identified with a structure $\cal C$ in $\cal K$, the induced topology $K_{L}|_{\cal C}$ can be identified with $J_{\cal C}$, the inclusion ${\cal C}\hookrightarrow L$ is a morphism of sites $({\cal C}, J_{\cal C})\to (L, K_{L})$ and we have a geometric morphism  
\[
\Sh(L, K_{L}) \to \Sh({\cal C}, J_{\cal C});
\] 
we can thus expect to be able to use the technique for generating reflections of \cite{OC13} to obtain a reflection or a coreflection from $\cal K$ to the resulting category of $\cal L$-structures which extends our duality.

As in the case of Stone and Priestley dualities, these categorical equivalences also admit purely topological formulations, arising from the fact that, by Theorem \ref{prgen}, the free structures under consideration admit topological representations as the $\cal L$-structures generated by the structures in $\cal K$ inside the powerset on the set of points of the associated toposes. Specifically, the free functors can be described topologically as the functors which assign to any structure $\cal C$ in $\cal K$ the $\cal L$-structure generated by inside the powerset on the set of points of the topos $\Sh({\cal C}, J_{\cal C})$, acting on the arrows accordingly.  

Concerning the problem of characterizing the extended images of these duality functors, we have the following result, which represents an analogue of Corollary \ref{corollarycar} in the context of the structures satisfying the hypotheses of the last part of Theorem \ref{prgen}. In order to state the theorem, we introduce some definitions. Given an invariant $C$ of families of subterminals in a topos satisfying the hypotheses of Theorem 3.25 \cite{OC11}, an element $l\in L$ is said to be $C$-compact if and only if it is $C$-compact as an open set of the space $Y_{L}$ of points of the topos $\Sh(L, K_{L})$ (via the canonical embedding $L\mono Id_{K_{L}}(L)\cong {\cal O}(Y_{L})\subseteq \mathscr{P}(Y_{L})$), equivalently if every $K_{L}$-covering sieve on $l$ admits a refinement in $L$ by a family satisfying $C$.

\begin{theorem}\label{free}
Let $\cal C$ be a preordered structure, $\cal L$ be a category of preordered structures and $\cal M$ be a class of functions from $\cal C$ to structures in $\cal L$. Let assume that all the powersets can be made into structures in $\cal L$ and that for any function $f:X\to Y$ the induced function $\mathscr{P}(f)=f^{-1}:\mathscr{P}(Y)\to \mathscr{P}(X)$ is an arrow in $\cal L$. Let us moreover assume that, given a function from $\cal C$ to a structure in $\cal L$, if the composition of it with each of the arrows belonging to a jointly injective family of arrows in $\cal L$ belongs to $\cal M$ then the function itself belongs to $\cal M$. 
Suppose that there exists a subcanonical Grothendieck topology $J_{\cal C}$ on $\cal C$ which is $C$-induced for an invariant $C$ satisfying the hypotheses of Theorem 3.25 \cite{OC11} and is such that the morphisms ${\cal C}\to \{0,1\}$ in $\cal M$ coincide precisely with a designated set $X_{\cal C}$ of jointly conservative flat $J_{\cal C}$-continuous functors ${\cal C}\to \{0,1\}$. Suppose moreover that the free $({\cal M}, {\cal L})$-structure on $\cal C$ exists and that the $\cal L$-substructure of $\mathscr{P}(X_{\cal C})$ generated by $\cal C$ (regarded as a subset of $\mathscr{P}(X_{\cal C})$ via the composite map ${\cal C}\mono Id_{J_{\cal C}}({\cal C})\cong {\cal O}(X_{\cal C})\subseteq \mathscr{P}(X_{\cal C})$) is contained in ${\cal O}(X_{\cal C})$. 

Let $i:{\cal C}\hookrightarrow L$ be an inclusion of $\cal C$ into a structure $L$ in $\cal L$. Suppose that $K_{L}$ is a subcanonical topology on $L$ such that the topos $\Sh(L, K_{L})$ has enough points. Then the following conditions are equivalent.
\begin{enumerate}[(i)]
\item $i$ realizes $L$ as the free $({\cal M}, {\cal L})$-structure on $\cal C$;

\item $i$ identifies $\cal C$ with the subset of elements of $L$ which are $C$-compact, $\cal C$ is $K_{L}$-dense in $L$, $J_{\cal C}=K_{L}|_{\cal C}$ and $L$ can be identified with the $\cal L$-structure generated by $\cal C$ inside $\mathscr{P}(Y_{L})$, where $Y_{L}$ is the space of points of the topos $\Sh(L, K_{L})$ and $\cal C$ is regarded as a subset of $\mathscr{P}(Y_{L})$ via the composite of $i$ with the canonical inclusion $L\mono {\cal O}(Y_{L})\subseteq \mathscr{P}(Y_{L})$.
\end{enumerate}
\end{theorem}

\begin{proofs}
$(i)\imp (ii)$. This follows from Theorem \ref{prgen}, specifically from the equivalence of toposes
\[
\Sh({\cal C}, J_{\cal C})\simeq \Sh(L, K_{L}).
\] 
$(ii)\imp (i)$. Under the hypotheses of condition $(ii)$, the Comparison Lemma ensures that the geometric morphism 
\[
\Sh(L, K_{L}) \to \Sh({\cal C}, J_{\cal C})
\]
induced by the morphism of sites $i:({\cal C}, J_{\cal C})\rightarrow (L, K_{L})$ is an equivalence of toposes. From this it follows that the morphisms $L\to \{0,1\}$ in $\cal L$ which extend the morphisms ${\cal C}\to \{0,1\}$ in $\cal M$ via $i$ coincide with a designated set $Y_{L}$ of jointly conservative flat $K_{L}$-continuous functors $L \to \{0,1\}$. Therefore the hypotheses of Theorem \ref{prgen} are satisfied and hence the $\cal L$-structure generated by $\cal C$, that is $L$, coincides with the free $({\cal M}, {\cal L})$-structure on $\cal C$.  
\end{proofs}

Before turning our attention to the examples, we remark that the construction of the free $({\cal L}, {\cal M}_{\cal C})$-structure $L$ on a structure $\cal C$, where $\cal L$ is a category of structures and ${\cal M}_{\cal C}$ is a class of maps from $\cal C$ to structures in $\cal L$ can, under appropriate conditions, be made functorial and give rise to a functor which is left adjoint to a forgetful one. Specifically, suppose that we have a category $\cal K$ of structures $\cal C$, and a forgetful functor $U:{\cal L}\to {\cal K}$ such that for any $L$ in ${\cal L}$ and any $\cal C$ in $\cal K$ the arrows ${\cal C} \to U(L)$ in $\cal K$ coincide exactly with the arrows ${\cal C}\to L$ in ${\cal M}_{\cal C}$. Then one can define a functor $F:{\cal K}\to {\cal L}$ which is left adjoint to $U:{\cal L}\to {\cal K}$, as follows. For any $\cal C$ in $\cal K$ we set $F({\cal C})$ equal to the free $({\cal L}, {\cal M}_{\cal C})$-structure on $\cal C$, while for any arrow $f:{\cal C}\to {\cal C}'$ in $\cal K$, we consider its composite with the canonical map ${\cal C}'\to U(L_{{\cal C}'})$ and set $F(f)$ equal to the unique arrow $g:L_{\cal C}\to L_{{\cal C}'}$ in $\cal L$ such that, denoted by $i_{\cal C}:{\cal C}\to U(L_{\cal C})$ the canonical map, $U(g)\circ i_{\cal C}=i_{{\cal C}'}\circ f$ (note that such an arrow exists and is unique by the universal property of the free $({\cal L}, {\cal M}_{\cal C})$-structure on $\cal C$). It is immediate to see that $F$ is left adjoint to $U$. Indeed, for any $\cal C$ in $\cal K$ and any $L$ in $\cal K$ the bijective correspondence between arrows ${\cal C}\to U(L)$ in $\cal K$ and arrows $F({\cal C})=L_{\cal C}\to L$ in $\cal L$ given by the universal property of the free $\cal L$-structures $L_{\cal C}$ is clearly natural in ${\cal C}\in {\cal K}$ (by definition of the functor $F$) and in $L\in {\cal L}$ by the universal property of the structures $L_{\cal C}$. Thus, under these hypotheses, we have a right adjoint to the free functor $F$ defined on an appropriate subcategory of $\cal L$, which yields an inverse to $F$ when restricted to the extended image $ExtIm(F)$ of $F$.

\section{Examples}\label{ex}

\subsection{A Priestley-type duality for coherent posets}\label{cohposets}

Let us recall the following notion from \cite{OC11}.

\begin{definition}
Let $\cal C$ and $\cal D$ be preorders. A monotone map $f:{\cal C}\to {\cal D}$ is said to be \emph{flat} if the following two conditions hold:
\begin{enumerate}[(i)]
\item For any $d\in {\cal D}$ there exists $c\in {\cal C}$ such that $d\leq f(c)$;

\item For any element $d\in {\cal D}$ and any elements $c,c'\in {\cal C}$ such that $d\leq f(c)$ and $d\leq f(c')$ there exists $c''\in {\cal C}$ such that $c''\leq c$, $c''\leq c'$ and $d\leq f(c'')$.
\end{enumerate}
\end{definition}

Let $\textbf{Pos}_{f}$ denote the category of posets and flat maps between them. By the results in \cite{OC11}, we have a duality between $\textbf{Pos}_{f}$ and a category $\cal V$ of topological spaces, namely the subcategory $\textbf{SSComp}$ of $\textbf{Top}$ whose objects are the sober spaces which a basis of supercompact open sets and whose arrows are the continuous maps between them such that the inverse image of any supercompact open set is supercompact. 

This duality sends a poset $\cal P$ in $\textbf{Pos}_{f}$ to the soberification of the Alexandrov space associated to its opposite, that is to the topological space $A({\cal P})$ whose underlying set if the set ${\cal F}^{dir}_{{\cal P}^{\textrm{op}}}$ of all the non-empty directed ideals on ${\cal P}^{\textrm{op}}$ and as open sets the subsets of the form
\[
{\cal F}_{U}=\{F\in {\cal F}^{dir}_{\cal P} \textrm{ | } F\cap U\neq \emptyset\},
\] 
where $U$ ranges among the lower sets in $\cal P$, and it sends a monotone map $f:{\cal P}\to {\cal Q}$ in $\textbf{Pos}_{f}$ to the map $A({\cal Q})={\cal F}^{dir}_{{\cal Q}} \to A({\cal P})={\cal F}^{dir}_{{\cal P}}$ sending a ideal set $I$ in ${\cal F}^{dir}_{{\cal P}^{\textrm{op}}}$ to the ideal in ${\cal F}^{dir}_{{\cal Q}^{\textrm{op}}}$ given by the inverse image $f^{-1}(I)$ of $I$ under $f$. Note that under this duality, finite posets correspond precisely to the finite and sober (equivalently, finite and $T_{0}$) topological spaces. 
 
Let $\cal P$ be a poset. It is shown in \cite{flatcoh} that the theory of flat functors on ${\cal P}$ (equivalently, the theory ${\mathbb T}_{\cal P}$ of non-empty directed ideals on ${\cal P}^{\textrm{op}}$ defined in \cite{OC11}) is coherent if and only if ${\cal P}$ has all fc finite limits (i.e., for any finite diagram $D$ in ${\cal P}^\textrm{op}$ there exists a finite family of cones on $D$ such that any other cone on $D$ factors through one in that family). We define a poset $\cal P$ to be \emph{coherent} if this condition is satisfied. Clearly, any finite poset is coherent. We denote by $\textbf{Pos}_{coh}$ the full subcategory of $\textbf{Pos}_{f}$ on the coherent posets. 

Recall from \cite{OC11} that the theory ${\mathbb T}_{\cal P}$ is defined over a signature consisting of one atomic proposition $F_{a}$ for each element $a\in {\cal P}$ and has the following axioms:
\[
(\top \vdash \mathbin{\mathop{\textrm{\huge $\vee$}}\limits_{c\in {\cal C}}} F_{c});
\]
\[
(F_{a} \vdash F_{b})
\] 
for any $a\leq b$ in ${\cal P}$; 
\[
(F_{a}\wedge F_{b} \vdash \mathbin{\mathop{\textrm{\huge $\vee$}}\limits_{c\in K_{a,b}}} F_{c})
\]
for any $a, b \in {\cal P}$, where $K_{a,b}$ is the collection of all the elements $c\in {\cal P}$ such that $c\leq a$ and $c\leq b$ in ${\cal P}$.

The fact that $\cal P$ is coherent means that all the sequents of the form 
\[
(\top \vdash \mathbin{\mathop{\textrm{\huge $\vee$}}\limits_{c\in {\cal C}}} F_{c});
\]
or of the form 
\[
(F_{a}\wedge F_{b} \vdash \mathbin{\mathop{\textrm{\huge $\vee$}}\limits_{c\in K_{a,b}}} F_{c})
\]
can be replaced by finitary (in fact, coherent) sequents over the same signature which are logically equivalent (in geometric logic) to them. Indeed, one can easily see that a poset $\cal P$ is coherent if and only if there exists a finite set $\{c_{1}, \ldots, c_{n}\}$ of elements of $\cal P$ such that for any $c\in {\cal P}$, $c\leq c_{i}$ in $\cal P$ for some $i\in \{1, \ldots, n\}$, and for any $c,c'\in {\cal P}$ there exists a finite set $\{d_{1}, \ldots, d_{m}\}$ of elements of $\cal P$ such that for any element $d\in {\cal P}$ such that $d\leq c$ and $d\leq c'$, $d\leq d_{j}$ in $\cal P$ for some $j\in \{1, \ldots, m\}$. 

For any coherent category $\cal D$ the models of the theory ${\mathbb T}_{\cal P}$ in $\cal D$ can be identified with the monotone maps $f:{\cal P}\to {\cal D}$ such that 
\begin{enumerate}[(i)]

\item for any (equivalently, every) finite set $\{c_{1}, \ldots, c_{n}\}$ in $\cal P$ such that for any $c\in {\cal P}$, $c\leq c_{i}$ in $\cal P$ for some $i\in \{1, \ldots, n\}$, $1=f(c_{1})\vee \cdots \vee f(c_{n})$ in $\cal D$, and 

\item for any $c,c' \in {\cal P}$ and any finite set $\{d_{1}, \ldots, d_{m}\}$ of elements of $\cal P$ such that for any element $d\in {\cal P}$ such that $d\leq c$ and $d\leq c'$, $d\leq d_{j}$ in $\cal P$ for some $j\in \{1, \ldots, m\}$, $f(c)\wedge f(c')=f(d_{1})\vee \cdots \vee f(d_{n})$.  
\end{enumerate}

If we define $\cal L$ to be the category of distributive lattices (resp. the category of Boolean algebras) and ${\cal M}_{\cal P}$ to be the class of morphisms from a given coherent poset $\cal P$ to structures of $\cal L$ which are models of ${\mathbb T}_{\cal P}$ then the coherent syntactic category $D_{\cal P}$ of the theory ${\mathbb T}_{\cal P}$ (resp. the coherent syntactic category $B_{\cal P}$ of the Morleyization of the theory ${\mathbb T}_{\cal P}$) can be identified, together with the canonical map $j_{\cal P}:{\cal P}\to D_{\cal P}$ (resp. $i_{\cal P}:{\cal P}\to B_{\cal P}$), with the free $({\cal L}, {\cal M}_{\cal P})$-structure on $\cal P$. In other words, $D_{\cal P}$ satisfies the universal property that for any map $f:{\cal P}\to {\cal D}$ in ${\cal M}_{\cal P}$ to a distributive lattice $\cal D$ there exists a unique distributive lattice homomorphism $\overline{f}:D_{\cal P}\to {\cal D}$ such that $\overline{f}\circ j_{\cal P}=f$; similarly, $B_{\cal P}$ satisfies the universal property that for any map $f:{\cal P}\to {\cal B}$ in ${\cal M}_{\cal P}$ to a Boolean algebra $\cal B$ there exists a unique homomorphism of Boolean algebras $\overline{f}:B_{\cal P}\to {\cal B}$ such that $\overline{f}\circ i_{\cal P}=f$.

The duality between $\textbf{Pos}_{f}$ and $\textbf{SSComp}$ clearly restricts to a duality between $\textbf{Pos}_{coh}$ and the full subcategory $\cal U$ of $\textbf{SSComp}$ on the topological spaces such that their collection of supercompact open sets forms, with the natural inclusion order, a coherent poset. 

By using the method of section \ref{genmeth}, we can lift this duality 
\[
A:{\textbf{Pos}_{coh}}^{\textrm{op}} \simeq {\cal U}
\]
to an duality 
\[
B:{\textbf{Pos}_{coh}}^{\textrm{op}}  \simeq {\cal V}
\]
between $\textbf{Pos}_{coh}$ and a category of preordered topological spaces $\cal V$.

Specifically, for any $\cal P$ in $\textbf{Pos}_{coh}$ we define $B({\cal P})$ as the preordered topological space whose underlying set is the set ${\cal F}^{dir}_{{\cal P}^{\textrm{op}}}$, endowed with the $R_{\cal P}$-patch topology, where $R_{\cal P}$ is the subset of $\mathscr{P}({\cal F}^{dir}_{{\cal P}^{\textrm{op}}})$ given by the composite ${\cal P}\mono Id({\cal P})\cong {\cal O}(A({\cal P}))\subseteq \mathscr{P}({\cal F}^{dir}_{{\cal P}^{\textrm{op}}})$, and with the specialization preorder corresponding to the space $A({\cal P})$. 

Given a Priestley space, we call a clopen upper set which cannot be decomposed as a proper union of clopen upper sets \emph{weakly indecomposable}. Since by Priestley duality the clopen upper sets form a basis for the coherent topology associated to the Priestley space, whose open sets are precisely the upper open sets, the upper (cl)open sets which are supercompact among the upper open sets coincide with the weakly indecomposable clopen upper sets.  

By the method of section \ref{genmeth}, the assignment ${\cal P}\rightarrow B({\cal P})$ can be made functorial and yields an equivalence between ${\textbf{Pos}_{coh}}^{\textrm{op}}$ and the subcategory of $\textbf{PTop}$ whose objects are the Priestley spaces $(X, \tau, \leq)$ such that for any $x, y\in X$ with $x\nleq y$, there is a clopen weakly indecomposable $\leq$-upper set $U$ of $\tau$ with the property that $x\in U$ and $y\notin U$ and whose arrows are the continuous order-preserving maps between them such that the inverse image of any weakly indecomposable upper clopen set is weakly indecomposable. 

Notice that under this equivalence, the finite posets correspond precisely to the finite sober ordered topological spaces $(X, \tau, \leq)$ such that for any $x, y\in X$ such that $x\nleq y$, there is a weakly indecomposable clopen $\leq$-upper set $U$ of $\tau$ with the property that $x\in U$ and $y\notin U$.       

\subsection{A Priestley-type duality for meet-semilattices}

We can obtain a Priestley-type duality for meet-semilattices by restricting the duality for coherent posets established above. Indeed, the category $\textbf{MsLat}$ of meet-semilattices and meet-semilattice homomorphisms between them can be identified as the full subcategory of the category ${\textbf{Pos}_{coh}}$ on the meet-semilattices. We thus obtain a duality
\[
B:\textbf{MsLat}^{\textrm{op}}\to \textbf{PTop}
\]
between $\textbf{MsLat}$ and a category of Priestley spaces, given by the extended image of the functor $B$.

Given a meet-semilattice $\cal M$, the Priestley space $B({\cal M})$ is the ordered topological space whose underlying set is the collection $X_{\cal M}$ of all the filters on $\cal M$, endowed with the topology generated by the sets of the form $\{F\in X_{\cal M} \textrm{ | } m\in F\}$ and their complements in $\mathscr{P}(X_{\cal M})$, and with the order $\leq_{\cal M}$ defined as follows: for any $F, F'\in X_{\cal M}$, $F\leq_{\cal M} F'$ if and only if $F\subseteq F'$. Given a meet-semilattice homomorphism $f:{\cal M}\to {\cal N}$, $B(f):X_{\cal N}\to X_{\cal M}$ is the map sending any filter $F$ in $X_{\cal N}$ to the filter in $\cal M$ given by the inverse image $f^{-1}(F)$. 

The extended image of the functor $B$ can be characterized as the category of ordered topological spaces whose objects are the Priestley spaces $(X, \tau, \leq)$ such that for any $x, y\in X$ with $x\nleq y$, there is a weakly indecomposable clopen $\leq$-upper set $U$ of $\tau$ with the property that $x\in U$ and $y\notin U$, and the intersection of any two weakly indecomposable clopen $\leq$-upper set is weakly indecomposable, and whose arrows are the continuous order-preserving maps between them such that the inverse image of any weakly indecomposable upper clopen set is weakly indecomposable.

The duality admits the following algebraic interpretation (cf. section \ref{algint} above).

We have a functor
\[
\tilde{B}:\textbf{MsLat} \to \textbf{Bool}_{\leq}
\]  

defined as follows: for any $\cal M$ in $\textbf{MsLat}$, $\tilde{B}({\cal M})=(B_{\cal M}, \leq_{\cal M})$, where $B_{\cal M}$ is the free Boolean algebra on $\cal M$ and $\leq_{\cal M}$ is the order on $Spec(B)$ given by: for any $F, F'\in Spec(B)$, $F\leq_{\cal M} F'$ if and only if $F\cap {\cal M}\subseteq F'\cap {\cal M}$, while for any arrow $f:{\cal M}\to {\cal M}'$ in $\textbf{MsLat}$, $\tilde{B}(f)$ is equal to the unique Boolean algebra homomorphism $B_{f}:B_{\cal M} \to B_{\cal N}$ which makes the following diagram commutes (where $i_{\cal M}:{\cal M}\to B_{\cal M}$ and $i_{\cal N}:{\cal N}\to B_{\cal N}$ are the canonical inclusions from the relevant meet-semilattice to the free Boolean algebra on it):

\[  
\xymatrix {
{\cal M} \ar[r]^{f} \ar[d]^{i_{\cal M}} & {\cal N} \ar[d]^{i_{\cal N}} \\
B_{\cal M} \ar[r]^{B_{f}} & B_{\cal V}}
\]

Any $\cal M$ in $\textbf{MsLat}$ can be recovered (up to isomorphism) from the Boolean algebra $B_{\cal M}$ as the set of its elements which are $\leq_{\cal M}$-upper and cannot be written as a proper finite join of $\leq_{\cal M}$-upper elements, while any arrow $f:{\cal M}\to {\cal N}$ in $\textbf{MsLat}$ can be recovered as the restriction of the corresponding arrow $B(f)$ to the subsets consisting of the upper elements which cannot be written as a proper finite join of upper elements. Given an object $(B, \leq)$ of $\textbf{Bool}_{\leq}$, we shall say that an element $b\in B$ is an indecomposable $\leq$-upper element if cannot be written as a proper finite join of upper elements in $B$.    

By Theorem \ref{tfae}, the extended image $ExtIm(\tilde{B})$ of the functor $\tilde{B}$ can be characterized as follows. The objects of $ExtIm(\tilde{B})$ are the objects $(B, \leq)$ of $\textbf{Bool}_{\leq}$ with the property that for any $F, F'\in Spec(B)$ such that $F\nleq F'$ there exists an indecomposable $\leq$-upper element $b\in B$ such that $b\in F$ and $b\notin F'$, and the meet of any two indecomposable $\leq$-upper elements is indecomposable; equivalently, the subset $B^{\ast}$ of $B$ consisting of the indecomposable $\leq$-upper elements $b$ of $B$, endowed with the induced order, is a meet-semilattice, the inclusion $B^{\ast}\subseteq B$ realizes $B$ as the free Boolean algebra on the meet-semilattice $B^{\ast}$ and for any $F, F'\in Spec(B)$ such that $F\nleq F'$ there exists an indecomposable $\leq$-upper element $b\in B$ such that $b\in F$ and $b\notin F'$. The arrows $(B, \leq)\to (B', \leq')$ in $ExtIm(\tilde{B})$ are the Boolean algebra homomorphisms $f:B\to B'$ such that $f^{-1}:Spec(B')\to Spec(B)$ is order-preserving and $f$ sends indecomposable $\leq$-upper elements of $B$ to indecomposable $\leq'$-upper elements of $B'$.

From Corollary \ref{corollarycar} we get the following criterion for a meet-semilattice inclusion $i:{\cal M}\to B$ into a Boolean algebra $B$ to realize $B$ as the free Boolean algebra on the meet-semilattice $\cal M$: $i$ realizes $B$ as the free Boolean algebra on $\cal M$ if and only if $\cal M$ can be identified, via $i$, with the subset of $B$ consisting of the indecomposable $\leq_{\cal M}$-upper elements of $B$.

\subsection{A Priestley-type duality for disjunctively distributive lattices}
 
Recall that a disjunctively distributive lattice is a meet-semilattice in which finite joins of pairwise disjoint elements exist and distribute over finite meets. 

In section 4 of \cite{OC11} we established a duality between the category $\textbf{DJLat}$ of disjunctively distributive lattices and meet-semilattice homomorphisms between them which send pairwise disjoint elements to pairwise disjoint elements and the category $\textbf{Top}_{dj}$ whose objects are the sober topological spaces with a basis of disjunctively compact open sets which is closed under finite intersections and satisfies the property that any covering of a basic open set has a disjunctively compact refinement by basic open sets, and whose arrows are the continuous maps between such spaces such that the inverse image of any disjunctively compact open set is a disjunctively compact open set. This duality sends a disjunctively distributive lattice $\cal D$ to the space of points of the topos $\Sh({\cal D}, J^{dj}_{\cal D})$, where $ J^{dj}_{\cal D}$ is the disjunctive topology on $\cal D$, and acts on the arrows in the natural way. 

We can lift this duality to a duality between $\textbf{DJLat}$ and a subcategory of $\textbf{PTop}$ by using the method of section \ref{genmeth}. Specifically, the Priestley space associated to a disjunctively distributive lattice $\cal D$ is the ordered topological space whose underlying set is the collection $X_{\cal D}$ of disjunctive filters on $\cal D$, endowed with the topology generated by the sets of the form $\{F\in X_{\cal D} \textrm{ | } d\in F\}$ and their complements in $\mathscr{P}(X_{\cal D})$ and with subset-inclusion ordering, while the map of Priestley spaces associated to a given homomorphism $f:{\cal D}\to {\cal D}'$ of disjunctively distributive lattices, is the map $X_{{\cal D}'}\to X_{\cal D}$ sending any filter $F$ in $X_{{\cal D}'}$ to the filter in $X_{{\cal D}}$ given by the inverse image $f^{-1}(F)$. 

A disjunctively distributive lattice can be recovered (up to isomorphism) from the associated Priestley space as the set of clopen uppper sets which satisfy the property that any covering of them by clopen upper sets admits a disjunctive refinement by clopen upper sets. Algebraically, a disjunctively distributive lattice $\cal D$ can be recovered from the free Boolean algebra $B_{\cal D}$ on it as the set of its $\leq$-upper elements (where $\leq$ is the order on the Priestley space associated to $\cal D$ under this duality) having the property that any finite covering of them by $\leq$-upper elements in $B_{\cal D}$ admits a disjunctive refinement by $\leq$-upper elements. 

From Corollary \ref{corollarycar}, we get the following criterion for a disjunctively distributive inclusion homomorphism $i:{\cal D}\to B$ into a Boolean algebra $B$ to realize $B$ as the free Boolean algebra on the disjunctively distributive lattice $\cal D$: $i$ realizes $B$ as the free Boolean algebra on $\cal D$ if and only if $\cal D$ can be identified, via $i$, with the subset of $B$ consisting of the $\leq_{\cal M}$-upper elements of $B$ with the property that any finite covering of them by $\leq$-upper elements in $B_{\cal D}$ admits a disjunctive refinement by $\leq$-upper elements.

\subsection{Other dualities}

\begin{enumerate}
\item \emph{Free distributive lattices on meet-semilattices}. 

By Theorem \ref{prgen}, for any meet-semilattice $\cal M$, denoted by $D_{\cal M}$ the free distributive lattice on $\cal M$, we have an equivalence of toposes
\[
\Sh(D_{\cal M}, J_{D_{\cal M}}) \simeq [{\cal M}^{\textrm{op}}, \Set],
\]
where $J_{D_{\cal M}}$ is the coherent topology on $D_{\cal M}$.

The meet-semilattice $\cal M$ can thus be recovered from $D_{\cal M}$ as the set of its elements which cannot be written as a proper finite join of elements of $D_{\cal M}$. Given a distributive lattice $\cal D$, we shall say that an element of $\cal D$ is indecomposable if it cannot be written as a proper finite join in $\cal D$ of elements of $\cal D$. 

The free functor $F:\textbf{MSLat}\to \textbf{DLat}$ thus has an inverse defined on its extended image, which sends a distributive lattice $D$ in $ExtIm(F)$ to the set of its indecomposable elements. This equivalence can be extended to a coreflection from $\textbf{MSLat}$ to the subcategory of $\textbf{DLat}$ whose objects are the distributive lattices whose set of indecomposable elements forms, with the induced order, a meet-semilattice, and whose arrows are the distributive lattice homomorphisms which send indecomposable elements to indecomposable elements. 

By Theorem \ref{free}, the extended image $ExtIm(F)$ of the functor $F$ can be characterized as the subcategory of $\textbf{DLat}$ whose objects are the distributive lattices $D$ such that any element can be written as a finite join of indecomposable elements and the meet of any two indecomposable elements is indecomposable and whose arrows are the distributive lattice homomorphisms between such lattices which send indecomposable elements to indecomposable elements. In fact, we have the following criterion for an meet-semilattice homomorphic inclusion $i:{\cal M}\hookrightarrow D$ of a meet-semilattice $\cal M$ into a distributive lattice to realize $D$ as the free distributive lattice on $\cal M$: $i$ realizes $\cal M$ as the set of indecomposable elements of $D$ and every element of $D$ can be written as a finite join of indecomposable elements.     
  
\item \emph{Free distributive lattices on disjunctively distributive lattices}. 

By Theorem \ref{prgen}, for any disjunctively distributive lattice $\cal D$, denoted by $D_{\cal D}$ the free distributive lattice on $\cal D$, we have an equivalence of toposes
\[
\Sh(D_{\cal D}, J_{D_{\cal D}}) \simeq \Sh({\cal D}, J^{dj}_{{\cal D}}),
\]
where $J_{D_{\cal M}}$ is the coherent topology on $D_{\cal M}$ and $J^{dj}_{\cal D}$ is the disjunctive topology on $\cal D$. 

The disjunctively distributive lattice $\cal D$ can thus be recovered from $D_{\cal D}$ as the set of its elements which are disjunctively compact, that is the elements $d$ of $D_{\cal D}$ such that for any finite family of elements $\{d_{1}, \ldots, d_{n}\}$ such that $d_{1} \vee \cdots \vee d_{n}=d$ there exists a finite family $\{a_{1}, \ldots, a_{m}\}$ of pairwise disjoint elements of $D_{\cal D}$ with the property that $a_{1} \vee \cdots \vee a_{m}=d$ and for any $i\in \{1, \ldots, m\}$ there exists an element $j\in \{1, \ldots, n\}$ such that $a_{i}\leq d_{j}$.    

The free functor $F:\textbf{DJLat}\to \textbf{DLat}$ thus has an inverse defined on its extended image, which sends a distributive lattice $D$ in $ExtIm(F)$ to the set of its disjunctively compact elements. This equivalence can be extended to a coreflection from $\textbf{DJLat}$ to the subcategory of $\textbf{DLat}$ whose objects are the distributive lattices whose set of disjunctively compact elements is closed under finite meets and whose arrows are the distributive lattice homomorphisms which send disjunctively compact elements to disjunctively compact elements. 

By Theorem \ref{free} and the results in \cite{OC11}, the extended image of the functor $F$ can be characterized as the subcategory of $\textbf{DLat}$ whose objects are the distributive lattices $\cal D$ such that there exists a set of elements $\cal B$ of $\cal D$ closed under finite meets with the property that every finite covering in $\cal D$ of element in $\cal B$ admits a disjunctive refinement by a family of elements of $B$ (i.e., for any finite family $\{e_{1}, \ldots, e_{m}\}$ of elements of $D$ such that $e_{1}\vee \cdots \vee e_{m}\in B$ there exists a subset of $B$ formed by pairwise disjoint elements each of which is less than or equal to $e_{j}$ for some $j\in \{1, \ldots, m\}$ and whose join is equal to $e_{1}\vee \cdots \vee e_{m}$), and whose arrows are the distributive lattice homomorphisms between such lattices which send disjunctively compact elements to disjunctively compact elements.

\item \emph{Free frames on posets}.

For any poset $\cal P$, we have an equivalence of toposes
\[
[{\cal P}^{\textrm{op}}, \Set]\simeq \Sh(A_{\cal P})
\]
where $A_{\cal P}$ is the frame of lower sets in $\cal P$ (cf. \cite{OC11}). From which it follows that $A_{\cal P}$ is the free $\cal M$-frame on $\cal P$ where $\cal M$ is the class of monotone maps from $\cal P$ to frames $F$ such that 

\begin{enumerate}[(i)]
\item $\mathbin{\mathop{\textrm{\huge $\vee$}}\limits_{p\in {\cal P}}}f(p)=1_{F}$; 

\item for any $a,b \in {\cal P}$, $f(a)\wedge f(b)=\mathbin{\mathop{\textrm{\huge $\vee$}}\limits_{p\in K_{a, b}}}f(p)$,
where $K_{a,b}$ is the collection of all the elements $c\in {\cal P}$ such that $c\leq a$ and $c\leq b$ in ${\cal P}$.  
\end{enumerate}

One can recover (up to isomorphism) any poset $\cal P$ from the associated frame $A_{\cal P}$ as the subsep of its supercompact elements, endowed with the induced order. 

In order to make this into a duality between posets and frames, one can proceed in two different ways. One possibility is to rely on the fact that any monotone map $f:{\cal P}\to {\cal Q}$ of posets $\cal P$ and $\cal Q$ induces a complete frame homomorphism (i.e., a frame homomorphism preserving arbitrary infima, that is having a right adjoint) $f^{-1}:A_{{\cal Q}}\to A_{\cal P}$ and can be recovered from it as the restriction of its right adjoint to the principal ideals. This gives rise to a duality between the category of posets and monotone maps between them and the category whose objects are the frames with a basis of supercompact elements and whose arrows are the complete homomorphisms between them (cf. section 4.1 of \cite{OC11} for the details). An alternative, covariant, way of functorializing the assignment ${\cal P}\to A_{\cal P}$ consists in considering the flat maps between posets and the frame homomorphisms which they induce. Specifically, any flat map $f:{\cal P}\to {\cal Q}$ between posets induces a frame homomorphism $A_{\cal P}\to A_{\cal Q}$, sending any lower set $I$ in $A_{\cal P}$ to the lower set in $A_{\cal Q}$ generated by the image $f(I)$ of $I$ under $f$ in $\cal Q$, from which the map $f$ can be recovered as its restriction to the subsets of principal ideals. This gives rise to an equivalence between the category of posets and flat maps between them and the category whose objects are the frames with a basis of supercompact elements and whose arrows are the frame homomorphisms between them which send supercompact elements to supercompact elements.
  
\end{enumerate}

\end{document}